\documentclass{paper}[12pt]

\usepackage{graphicx,tikz,tikz-cd,psfrag}
\usetikzlibrary{matrix,arrows,decorations.pathmorphing}
\usepackage{amssymb,amsmath,mathtools}

\usepackage{setspace}

\usepackage{mathpazo}

\newtheorem{theorem}{\bf Theorem}
\newtheorem{lemma}{\bf Lemma}

\newenvironment{itemequation}{\vspace*{-\baselineskip}\[}{\]}

\newenvironment{itemeqnarray}{\vspace*{-\baselineskip}\eqnarray}{\endeqnarray}

\newenvironment{remark}{{\bf Remark }}

\newenvironment{proof}{{\bf Proof: }}

\newenvironment{definition}{{\bf Definition }}

\newenvironment{corollary}{{\bf Corollary }}

\newcommand{\modi}{\mbox{ $(\textrm{mod } \mathcal I)$}}

\title{A quotient framework for single-input pole placement and associated algorithms}

\author{Ph. M\"ullhaupt}

\begin{document}

\maketitle

\tableofcontents

\begin{abstract}
Single-input eigenvalue assignement problem (SEVAS) for dense (non-sparse) weakly controllable pairs $A$, $B$ using Ackermann's formula is revisted. Factorizations are presented that 
are interpreted using quotient (factor) vector spaces. Depending on the representations chosen for the equivalence class (given as specific projections by real orthogonal matrices), the numerical behavior of the pole placement scheme can be enhanced.  One version operates by placing one pole at a time (with complex conjugate poles grouped together to avoid complex arithmetic). Another version  operates with the coefficients of the characteristic polynomial directly. The latter version uses orthogonal real matrices and does not require complex-numbers floating-point arithmetic. In both cases, there are no constraint on the number of identical poles. A version uses only ring arithmetic. The algorithms are compared with numerically stable algorithms that appeared in the literature, such as the Miminis-Paige algorithm and the Varga pole shifting method. The justification rests on computer comparision using floating-point arithmetic of different precisions.
\end{abstract}

\section{Introduction}

It was in 1988 that George Miminis and Chris Paige published  a multi-input numerically stable algorithm eigenvalue assignment, which was shown, on very specific problems, to outperform the robust pole placement of Kautsky-Nichols-Van Dooren published in 1984 \cite{KautskyNicholsVanDooren}. The latter is still the main method in commercial numerical software for eigenvalue assignement (e.g. Matlab R2023b), and works very well for general systems. The main objective and motivation is to ensure some robustness to the initial data. However, in the single-input case, where the exact solution is uniquely defined, the Kautsky-Nichols-Van Dooren method could not improve on the two existing numerically stable algorithms available at that time, namely the pole shifting technique of Varga \cite{Varga}  based on the real Schur form and the single-input technique of Miminis-Paige \cite{MiminisPaige1} based on a preliminary Hessenberg reduction.

Other methods for single-input eigenvalue assignment (and also multi-input) exist such as Arnold and Datta \cite{ArnoldDatta2}, Bhattacharyya and DeSouza \cite{BhattacharyyaDeSouza}, Bru, Mas, and Urbano \cite{Urbano}, Bru, Cerdan, and Urbano \cite{Urbano2}, Datta \cite{Datta}, Patel and Misra \cite{PatelMisra}, Petkov, Christov, and Konstantinov \cite{Petkov}-\cite{Petkov2}, Tsui \cite{Tsui}. All algorithms that are numerically stable are compared from the perspective of QR decomposition by Arnold and Datta in \cite{ArnoldDatta} using a unified RQ reformulation.

In this present work, we develop an understanding of the pole placement for single-input linear time-invariant systems form the perspective of equivalent classes and quotients, taking advantage of the geometry, and not necessarily relying on orthogonal tranformations. This allows revisiting Ackermann's formula from the perspective of orthogonal transformations and non-orthogonal transformations, connecting the theory to a Lie-Algebroid \cite{Mullhaupt1}-\cite{Mullhaupt2} used to describe feedback linearization of affine single-input nonlinear systems \cite{Mullhaupt3}. We obtain a geometric understanding of the vector of gains (when the eigenvalues are real)  as the intersection of affine hyperplanes that are not necessarily orthogonal, each plane orientation depending on the real eigenvalue to be placed and the row geometry of the $A$ matrix and the orientation of the $B$ vector.

Numerical conditionning is explored through simulation and examples. The state-space  structure of the system is known to machine precision (matrices $A$ and $B$) and can undergo a similarity transform of the orthogonal type to within machine precision. This stage setting is quite different from the robustness criteria that the algorithm presented in \cite{KautskyNicholsVanDooren} had to withstand. The aim is to ensure closed-loop stability in simulation based on exact knowledge of the structure given within machine precision.

\part{Theoretical material and  results}

\section{Quotients, algebroids, projections (anchors)}

We refer to \cite{Gantmacher1}, p. 175, Chapter VII, for the definition of quotients and factor vector spaces.

\subsection{Congruences and factor vector spaces}

Consider an abstract vector space $\mathcal R$ together a vector subspace $\mathcal I \subset \mathcal R$. 
\begin{definition}
Two vectors $x \in \mathcal R$ and $y \in \mathcal R$ are congruent modulo $\mathcal I$ written $x \equiv y \mbox{ $(\textrm{mod } \mathcal I)$}$ if and only if $x - y \in \mathcal I$. 
\end{definition}

This is an equivalence relation since it has the properties of reflexivity, symmetry and transitivity, namely
\begin{enumerate}
\item $x \equiv x \modi$
\item $x \equiv y \modi \Rightarrow y \equiv x \modi$
\item $x \equiv z \modi  \textrm{ and } z \equiv y \modi \Rightarrow x \equiv z \modi$
\end{enumerate}

To any $x$ one can associate the set
\[ \bar x = 
\{ v | v \equiv x \modi \}
\]
Notice that $\bar x$ is not defined through a unique $x$ but any vector $v \equiv x \modi$ can stand for a representative of $\bar x$ and it is not difficult to show that indeed $\bar x$ does not depend on the representative $x$ chosen.
Such sets are called equivalence classes (shortly classes) and denoted as $\{\bar x \}$.

\begin{theorem}
In the manifold $\bar {\mathcal V}$ of all classes $\bar x$, $\bar y$, $\ldots$ two operations are introduced $\alpha \bar x$ and $\bar x + \bar y$ that turn these classes into a vector space.
\end{theorem}

This is a classical result that state that each class constitutes a specific vector in a new vector space, the factor (quotient) vector space written $\mathcal R / \mathcal I$. The operations are given as
\begin{eqnarray}
\alpha \bar x &=& \{ v | v \equiv \alpha x \modi \} \nonumber \\
\bar x + \bar y &=& \{ v | v \equiv x + y \modi \} \nonumber
\end{eqnarray}
and they do not depend on the representatives $x$ and $y$ chosen to denote the class $\bar x$ and $\bar y$.

\subsubsection{Abstract operators}

Let us associate an abstract operator $\mathfrak A$ to the abstract vector space $\mathcal R$ and let us pick a vector $B \in {\mathcal R}$ and let us set ${\mathcal I} := {\mathcal B} = \mbox{ span } B = \{ v | v = \lambda B, \forall \lambda \in \mathbb{R} \}$

\begin{tikzcd}
\mathcal R \arrow{r}{\mathfrak A} \arrow{rd} 
& \mathcal R \arrow{d}{\pi} \\
& {\mathcal R}/ {\mathcal B}
\end{tikzcd}

The homomorphism $\pi$ is the canonical homomorphism that associates to any vector $x \in \mathcal R$ its associated equivalence class $\bar x$ in $\mathcal R / \mathcal B$. 

Now, let us consider $\mathfrak A B \in \mathcal R$ and set 
\[ \mathcal I = \mbox{ span }\{B, \mathfrak A B\} = \{ v | v = \lambda_1 B + \lambda_2 \mathfrak A B, \forall \lambda_1, \lambda_2 \in \mathbb{R} \} \]

\begin{tikzcd}
\mathcal R \arrow{r}{\mathfrak A} \arrow{rd}{\Phi_1} & \mathcal R \arrow{r}{\mathfrak A} \arrow{d}{\pi_1} \arrow{r}{\mathfrak A} & \mathcal R \arrow{d}{\pi_1} \arrow{r}{\mathfrak A} 
& \mathcal R \arrow{d}{\pi_1} \\
 & \mathcal R/ \mbox{ span }\{B\} \arrow{rd}{\Phi_2} & \mathcal R/ \mbox{ span }\{B \} \arrow{d}{\pi_2}  & \mathcal R/ \mbox{ span }\{B\} \arrow{d}{\pi_2}  \\
& & \mathcal{R}/ \mbox{ span }\{B, {\mathfrak A}B\} \arrow{rd}{\Phi_3} & \mathcal{R}/ \mbox{ span }\{B, {\mathfrak A}B\} \arrow{d}{\pi_3} \\
  &  & & \mathcal{R}/ \mbox{ span }\{B, {\mathfrak A}B, {\mathfrak A}^2B\} 
\end{tikzcd}

%\begin{tikzcd}
%\mathcal R \arrow{r}{\mathfrak A} \arrow{rd}
%&
%\mathcal R \arrow{r}{\mathfrak A} %\arrow{rd}%{\Phi_1} %&  {\mathcal R}  \arrow{r}{\mathfrak A} \arrow{d}{\Pi_1} & 
%{\mathcal R} \arrow{r}{\mathfrak A} \arrow{d}{\Pi_1} \\
%& {\mathcal R}/ \mbox{ span }\{B\} \\arrow{rd}{\Phi_2} & \arrow{d}{\pi_2} & {\mathcal R}/\mbox{ span }\{B\} \arrow{d}{\pi_2} 
%\end{tikzcd}

%%%%%%%%%%%%%%%%%%%%%%%%%%%%%%%%
% some diagrams for testing
%
%\begin{tikzcd}
%A \arrow{r}{a} \arrow{d}{b}
%&B \arrow{d}{c}\\
%C \arrow{r}{d} &D
%\end{tikzcd}
%
%\begin{tikzcd}
%A \arrow[hook]{r}{u}[swap]{b}
%\arrow[two heads]{rd}{u}[swap]{b}
%&B \arrow[dotted]{d}{r}[swap]{l}
%\arrow[hookleftarrow]{r}{u}[swap]{b}
%&c \Arrow[two heads]{ld}{b}[swap]{u}\\
%&D
%%\end{tikzcd}
%
%\begin{tikzcd}
%A \arrow[hook]{r}\arrow[two heads]{rd}
%&B \arrow[dotted]{d}\arrow[hookleftarrow]{r}
%&C \arrow[two heads]{ld}\\
%&D
%\end{tikzcd}
%
%%\begin{tikzcd}
%&A\arrow{ldd}[swap]{f}\arrow{rd}[description]{c}
%\arrow{rrd}[description]{d}
%\arrow{rrrd}[description]{e}\\
%&B\arrow{ld}\arrow{r}&C\arrow{r}&D\arrow{r}&E\\
%F
%\end{tikzcd}

\subsubsection{Representations of a factor vector space in the ambient vector space}

\label{AppendixRepresentation}
Let us assume that the dimension of $\mathcal R$ is $n$.
A scalar product $<.|.>$ is chosen which is a bilinear map
\[
<.|.> : {\mathcal R} \times {\mathcal R} \rightarrow \mathbb{R}
\]
that satisfies the Schwarz inequality
\[
<x|y> \leq <x|x> + <y|y>
\]
Once this scalar product is chosen, a maximal set of independent vectors $e_1$, $e_2$, $\ldots$, $e_n$, all belonging to $\mathcal R$, can be chosen such that 
\begin{eqnarray}
<e_i|e_j> &=& 0 \qquad i \neq j \nonumber \\
<e_i|e_i> &=1& \qquad i = j \nonumber 
\end{eqnarray}

The set of vectors $\{ e_1, \ldots, e_n \}$ with the above property is a maximal set of orthonormal vectors. Once this set of orthonormal vectors is chosen, 
the scalar product can be computed using only operations in $\mathbb R$

\begin{lemma}
Take any two vectors $v_1$ and $v_2$ of $\mathcal R$ and express them as linear combinations of the orthonormal basis $\{ e_1, e_2, \ldots, e_n \}$ as 
\begin{eqnarray}
v_1 &=& \sum_{i=1}^n \alpha_i e_i \nonumber \\
v_2 &=& \sum_{i=1}^n \beta_i e_i \nonumber
\end{eqnarray}
with $\alpha_i \in \mathbb{R}$ and $\beta_i \in \mathbb{R}$, $i=1,\ldots, n$,
then the scalar product of $v_1$ and $v_2$ can be expressed as
\[
<v_1|v_2> = \sum_{i=1}^n \alpha_i \beta_i
\]
\end{lemma}

The proof follows directly from the bilinearity property of the scalar product $<.|.>$ and the fact that the family $\{e_1, e_2, \ldots, e_n\}$ is orthonormal.

It is also possible to represent the linear map $\mathfrak A$ by a matrix, the columns of which are the images of each of the vectors of the orthonormal basis, once these images are expressed as a linear combination of the vectors of the basis
\[
{\mathfrak A} e_j = \sum_{i=1}^n <e_i|\mathfrak A e_j> e_i
\]
The numbers $a_{ij} = <e_i | \mathfrak A e_j> \in \mathbb{R}$ are the components of the matrix
\begin{eqnarray}
A &=& \left( \begin{array}{cccc} 
<e_1|{\mathfrak A} e_1 > & <e_1|{\mathfrak A} e_2 > &  \cdots & <e_1|{\mathfrak A} e_n> \\
<e_2|{\mathfrak A} e_1 > & <e_2|{\mathfrak A} e_2 > & \cdots & <e_2 | {\mathfrak A} e_n> \\
\vdots  & \vdots & \vdots  & \vdots \\
<e_n|{\mathfrak A} e_1 > & <e_n|{\mathfrak A} e_2 >  & \cdots & <e_n | {\mathfrak A} e_n> 
\end{array} \right) \nonumber 
\end{eqnarray}

Consider the diagram

\begin{tikzcd}
\mathcal R \arrow{r}{\mathfrak A} \arrow{rd}{\Phi} & \mathcal R \arrow{d}{\pi}\\
& \mathcal R/\mbox{ span }\{B\} 
\end{tikzcd}

Notice that once $B$ is represented by $n$ real numbers $\beta_1$, $\beta_2$, $\ldots$, $\beta_n$ through the choice of the orthonormal family $\{e_1, e_2, \ldots, e_n \}$, that is,
\[
B = \sum_{i=1}^n \beta_i e_i
\] 
it is also possible to find suitable representations of $\mathcal R / \mbox{ span } \{ B \}$.
We will now show that this is possible by choosing a set of $n-1$ linearly independent vectors $g_1$, $g_2$, $\ldots$, $g_{n-1}$ all belonging to $\mathcal R$.

\begin{lemma}
Any set of linearly independent vectors $g_1$, $g_2$, $\ldots$, $g_{n-1}$ belonging to the vector space $\mathcal R$ such that $g_1, g_2, \ldots, g_{n-1}, B$ is a basis of $\mathcal R$
can serve as a basis of a representation of the vector space $\mathcal R/ \mbox{ span } \{B\}$. In that case, the vector
\[
\bar x = \{ v | v = x + \lambda B, \lambda \in \mathbb{R} \} 
\]
is represented as $\sum_{i=1}^{n-1} \gamma_i g_i$ through $n-1$ numbers $\gamma_1$, $\gamma_2$, $\ldots$, $\gamma_{n-1}$ that do not depend on the representative $x$ chosen for $\bar x$.
\end{lemma}

\begin{proof}
We must explicitly construct the isomorphism between any class
\[
\bar x = \{ v | v = x + \lambda B, \lambda \in \mathbb{R} \}
\]
and the corresponding element $\sum_{i=1}^{n-1} \gamma_i g_i$ associated to $\bar x$ in a unique way.

Apart from $x$, the set $\bar x$ might as well have different representatives, say $x_2$, $x_3$, etc. These representatives will all differ from $x$ by a suitable multiple of $B$, that is $x_2 = x + \lambda_2 B$, $x_3 = x + \lambda_3 B$, etc. As a consequence, if $x$ is written with respect to the basis $\{g_1, \ldots, g_n, B \}$ as
\[
x = \sum_{i=1}^{n-1} \gamma_i g_i + \beta B
\]
then
\begin{eqnarray}
x_2 &=& \sum_{i=1}^{n-1} \gamma_i g_i + \beta B + \lambda_2 B \nonumber \\
x_3 &=& \sum_{i=1}^{n-1} \gamma_i g_i + \beta B + \lambda_3 B \nonumber
\end{eqnarray}
with the same coefficients $\gamma_1$, $\ldots$, $\gamma_{n-1}$. Therefore the class $\bar x$ is associated with a unique set of numbers $\gamma_1$, $\ldots$, $\gamma_{n-1}$.

Reciprocally, given $\gamma_1$, $\ldots$, $\gamma_{n-1}$, let us now show that there exists a single equivalence class $\bar x$ associated to it. The proof proceeds by contradiction.
Therefore, let us suppose that two distinct equivalence classes 
\[
\bar x = \{ v | v = x + \lambda B, \lambda \in \mathbb{R} \}
\]
and
\[
\bar y = \{ v | v = y + \lambda B, \lambda \in \mathbb{R} \}
\]
 with $x$ and $y$ as representatives, such that there is no $\mu \in \mathbb{R}$ for which $x = y + \mu B$.
Otherwise stated, the only way to write $0$ as a linear combination of $x$, $y$ and $B$ is with zero coefficients, i.e. 
\[
\alpha_1 x + \alpha_2 y + \alpha_3 B = 0 \quad \Rightarrow \quad \alpha_1=\alpha_2=\alpha_3=0
\]
Recall (first part of the proof of the forward statement) that to any equivalence class there is a single set of real numbers $\gamma_1$, $\ldots$, $\gamma_n$ associated to it, i.e. to $\bar x$ we have $\gamma_{1,1},\ldots,\gamma_{1,n-1}$ and to $\bar y$  we have $\gamma_{2,1}, \ldots, \gamma_{2,n-1}$. By hypothesis (initiating our contradictive argument) $\gamma_{1,i} = \gamma_{2,i}$, $i=1,\ldots, n-1$ which means that
\begin{eqnarray}
x &=& \sum_{i=1}^{n-1} \gamma_i g_i + \beta_1 B \nonumber \\
y &=& \sum_{i=1}^{n-1} \gamma_i g_i + \beta_2 B \nonumber
\end{eqnarray}
But these last two equalities imply that 
\[
y - x - (\beta_1 - \beta_2) B = 0
\]
leading to a contradiction. 
\end{proof}

\subsubsection{Geometric interpretation}

\begin{lemma}
The family of linearly independent vectors $\{g_1, \ldots, g_{n-1} \}$ lie in a hyperplane whose orthogonal is never orthogonal to the vector $B$, that is, 
\[
<n|g_i> = 0 \qquad i=1,\ldots, n-1
\]
with
\[
<n|B> \neq 0
\]
\end{lemma}

\begin{proof}
A direct consequence of the fact that $g_1$, $\ldots$, $g_{n-1}$ are linearly independent and that $<g_i|B> \neq 0$, $i=1, \ldots, n-1$.
\end{proof} 

\begin{corollary}
All representations of $\mathcal B/ \mbox{ span } B$ in $\mathbb{R}^n$ are determined by the determination of a single vector $n$. 
\end{corollary}

We therefore distinguish two classes of representations for $\mathcal R/ \mbox{ span } B$
\begin{enumerate}
\item Orthogonal representations, those for which $n$ is chosen such that $n = \mu B$ for a certain value $\mu \in \mathbb{R}$. 
\item Oblique representations, those for which $<n|B> \neq 0$ and $\eta n + \beta B = 0 \Rightarrow \eta = \beta = 0$.
\end{enumerate}

\subsection{Algebroids, commutators and anchors}

\label{anchors}
Once a representation is chosen for the factor space (see Section \ref{AppendixRepresentation}), one can compare the effect of the induced matrix in that representation and whether or not the commutativity of the associated diagram occurs. We show that the initial $A$ matrix has to be modified in order for the representation diagram to commute. This gives rise to a specific commutator operation inducing a Lie-algebroid bracket associated with the Lie-algebra bracket given as the commutator of matrices.

Given two vectorfields $m_1$ and $m_2$, a 1-form $\omega$ and a vectorfield $g$, define the following bracket

\begin{eqnarray}
\label{algebroidBrack}
<m_1,m_2> = [m_1,m_2] + {\omega m_2 \over \omega g} [g,m_1] - {\omega m_1 \over \omega g} [g,m_2]
\end{eqnarray}

It has been shown in \cite{Mullhaupt1, Mullhaupt2} that this bracket is a Lie-algebroid bracket with anchor map

\[
\mbox{an}_{g,\omega} = \mbox{Pr } \Phi^*
\]

\subsubsection{Orthogonal commutators}
In our setting the vectorfields $m_1 = A_1 x$ and $m_2 = A_2 x$
and the projection map is
\[
\mbox{Pr }\Phi^* m = Q^T A Q x
\]
where $Q^T$ is a $n-1 \times n$ matrix containing orthornormal vectors with $q^T$ the complementary orthormal vector so that 
\[
P = \left( \begin{array}{c} q^T \\ Q^T \end{array} \right)
\] 
is full rank
\[
P^T P = I
\]
We can then consider $\omega = q^T$ which gives the commutator of vectorfields

\[
<A_1 x, A_2 x> = [A_1 x, A_2 x] + (q^T A_1 x) A_2 q - (q^T A_2 x) A_1 q
\]

or in matrix form, through the definition of a new commutator $<A_1,A_2>$ between matrices
$A_1$ and $A_2$ (not vectorfields) defined as
\[
<A_1,A_2> = A_1 A_2 - A_2 A_1 + A_2 q q^T A_1 - A_1 q q^T A_2
\]
The following commutation relations hold
\[
Q^T <A_1,A_2> Q = [Q^T A_1 Q, Q^T A_2 Q]
\]
or in full
\[
Q^T (A_1 A_2 - A_2 A_1 + A_2 q q^T A_1 - A_1 q q^T A_2) Q = Q^T A_1 Q Q^T A_2 Q - Q^T A_2 Q Q^T A_1 Q
\]
Noticing that $q^T Q = 0$ by definition we can put the previous relation into the following form
\begin{eqnarray}
Q^T (A_1 (I - q q^T) A_2 (I - q q^T) - A_2 (I - q q^T) A_1 (I - q q^T)) Q = \nonumber \\
Q^T A_1 Q Q^T A_2 Q - Q^T A_2 Q Q^T A_1 Q \nonumber
\end{eqnarray}
which defines another new commutator $<<A_1,A_2>>$ defined as
\[
<<A_1,A_2>> = [A_1 (I - q q^T), A_2 (I - q q^T)]
\]
This means that we can have a  set of (orthogonal based) commutators parameterized by two real parameters $\alpha$ and $\beta$ as
\[
<<A_1,A_2>>_{\alpha,\beta} = [A_1 (I - q q^T) + \alpha q q^T, A_2 (I - q q^T) + \beta q q^T]
\]

A particular commutator in this last family of commutators
will be most helpful in analysing pole assignment algorithms.

\subsubsection{Example illustrating the commutator $<<A_1,A_2>>$}

A small script will illustrate the commutator relations. Because it resorts
to a qr decomposition, the result may depend on the implementation of the qr algorithm 
(i.e. using MATLAB instead of SysQuake the results differ due to the nonunicity of the qr
decomposition, while validating the commutator relation nevertheless).

\begin{verbatim}
A1 = [1 3 5; 7 13 17; 1 1 1];
A2 = [2 4 6; 13 3 1; 7 5 3];

[QQ,r] = qr([-21 -5 5; 1 38 49; -4 12 3]);

qT = QQ(1,:);
QT = QQ(2:3,:);

% let us check the commutator relation with

A1r = QT*A1*QT'
A2r = QT*A2*QT'

% Expression (A)
A1r*A2r-A2r*A1r

% Expression (B)
QT*(A1*A2-A2*A1 + A2*qT'*qT*A1 - A1*qT'*qT*A2)*QT'
\end{verbatim}

If expression (A) is identical to Expression (B), then the commutator relation holds.
MATLAB returns the following:

\begin{verbatim}
>> % Expression (A)
>> A1r*A2r-A2r*A1r

ans =

   16.7141   89.8467
   83.5973  -16.7141

>> % Expression (B)
>> QT*(A1*A2-A2*A1 + A2*qT'*qT*A1 - A1*qT'*qT*A2)*QT'

ans =

   16.7141   89.8467
   83.5973  -16.7141
\end{verbatim}

\subsubsection{Oblique commutators}

The anchor 
\[
\mbox{an}_{\omega,g} (m) = m -  \frac{\omega m}{\omega g} g
\]
where
$\omega$ is a 1-form and $g$ a vectorfield 
replaces the previous anchor. The bracket (\ref{algebroidBrack}) remains a Lie algebroid with this newly defined anchor \cite{Mullhaupt2}.
Let us translate this is linear vectorfields (matrix notation). Here $\omega^T$ is simply a row constant vector and $g$ a column constant vector, $A$ is a matrix.

\[ 
\mbox{an}_{\omega,g}(A) = A - g \frac{\omega^T A}{\omega^T g} 
\]

which is a matrix of same dimension as the $A$ matrix.

The oblique commutator between two matrices, parameterized by $\omega^T$ and $g$, is defined as
\[
\{ \{A_1, A_2 \} \} = A_1 A_2 - A_2 A_1 + \frac{1}{\omega^T g} A_2 g \omega^T A_1 - \frac{1}{\omega^T g} A_1 g \omega^T A_2
\]
Setting 
\[
G = \frac{g \omega^T}{\omega^T g}
\]
the previous expression writes
\[
\{ \{A_1, A_2 \} \} = A_1 A_2 - A_2 A_1 + A_2 G A_1 - A_1 G A_2
\]
with anchor
\[
\mbox{an}_{\omega,g}(X) = X - G X
\]

The algebroid property becomes 
\[
\mbox{an}_{\omega,g}(\{ \{A_1, A_2 \} \}) = [\mbox{an}_{\omega,g}(A_1), \mbox{an}_{\omega,g}(A_2)] 
\]
where $[A_1,A_2] = A_1 A_2 - A_2 A_1$ is the classical commutator of matrices.
The algebroid property is used in oblique pole placement.

\subsubsection{Example illustrating the commutator $\{\{A_1, A_2\}\}$}

Let 
\[
A_1 = \left( \begin{array}{ccc} 1 & 3 & 5 \\ 7 & 13 & 17 \\ 1 & 1 & 1 \end{array} \right) \qquad
A_2 = \left( \begin{array}{ccc} 2 & 4 & 6 \\ 13 & 3 & 1 \\ 7 & 5 & 3 \end{array} \right)
\]
and
\[
\omega^T = \left( \begin{array}{ccc} 1 & 2 & 3 \end{array} \right) \qquad g = \left( \begin{array}{c}
14 \\ -2 \\ -3 \end{array} \right)
\]
we have with $\mbox{an}_{\omega,g} A =  A - g \frac{\omega^T A}{\omega^T g}$
\[
\mbox{an}_{\omega,g} A_1 = \left( \begin{array}{ccc} -251 & -445 & -583 \\ 43 & 77 & 101 \\
55 & 97 & 127 \end{array} \right) \qquad
\mbox{an}_{\omega,g} A_2 = \left( \begin{array}{ccc} -684 & -346 & -232 \\ 111 & 53 & 35 \\
154 & 80 & 54 \end{array} \right)
\]
Considering 
\[G = \frac{1}{\omega^T g} \, g \omega^T =  \left( \begin{array}{ccc}  14 &  28 & 42 \\
    -2  &  -4 & -6 \\
    -3  &  -6 & -9 \\
\end{array} \right)
\]
we cross-check the commutator relation
\begin{eqnarray}
\mbox{an}_{\omega,g}(A_1 A_2 - A_2 A_1 + A_2 G A_1- A_1 G A_2 ) &=& \nonumber \\
 \mbox{an}_{\omega,g}(A_1) \mbox{an}_{\omega,g}(A_2)- \mbox{an}_{\omega,g}(A_2) \mbox{an}_{\omega,g}(A_1) 
&=& \nonumber \\ \left( \begin{array}{ccc} 
 -111539 & -238613 & -323187 \\
 18346 & 39202 & 53088 \\
 24949 & 53403 & 72337 \end{array} \right)  &&\nonumber
\end{eqnarray}

\section{Ackermann's formula}

Consider a controllable linear system of dimension $n$

\[
\dot x = A x + B u
\]

with a single input. We have the following well known formula \cite{Ackermann}

\begin{theorem}
The feedback gain $K$ of the control law $u = -K x$ is  given by
\begin{eqnarray}
K = e_n^T {\mathcal C}^{-1} \Phi(A)
\end{eqnarray}
with $e_n$ the canonical basis vector is such that 
\[
\left| A - B K - \lambda I \right| = \Phi(\lambda) 
\]
\end{theorem}

This is known as Ackermann's formula in the literature.

\subsection{Factorization of Ackermann's formula}

\label{AckermannFactorization}

The first factorization is described as a Theorem.

\begin{theorem}
Let $\lambda_1$, $\lambda_2$, $\ldots$, $\lambda_n$ be $n$ eigenvalues to be assigned by a properly chosen row vector of gains $K$.
Let $C$ denote the last line of the inverse of the controllability matrix
\[
C = e_n^T {\mathcal C}^{-1}
\]
where
\[
{\mathcal C} = \left( \begin{array}{cccc}
B & AB & \cdots & A^{n-1} B 
\end{array} \right)
\]
This means that when $C$ is such that
\[
C \mathcal C = \left( \begin{array}{ccccc}
0 & 0 & \cdots & 0 & 1 \end{array} \right)
\]
then the following recursive procedure (nested iterations) computes a suitable row vector of feedback gains
\begin{eqnarray}
K_0 &=& C  \label{induction1}  \\
K_1 &=& K_0 A - \lambda_1 K_0 \nonumber \\
K_2 &=& K_1 A - \lambda_2 K_1 \nonumber \\
\vdots && \vdots \nonumber \\
K_{n} &=& K_{n-1} A - \lambda_n K_{n-1} \nonumber
\end{eqnarray}
The gain is then given by $K := K_n$.
\end{theorem}

The proof rests partly on the following lemma

\begin{lemma}
Let the characteristic matrix polynomial be factorized as
\[
\Phi(A) = (A - \lambda_n I) (A - \lambda_{n-1} I) \cdots (A - \lambda_1 I)
\]

then $\Phi(A) = \Phi_n(A)$ can be obtained using the recursive Horner-like scheme 
\begin{eqnarray}
\Phi_0 &=& I \nonumber \\
\Phi_1 &=& A \Phi_0  - \lambda_1 \Phi_0 \nonumber \\
\Phi_2 &=& A \Phi_1 - \lambda_2 \Phi_1 \nonumber \\
\vdots && \vdots \nonumber \\
\Phi_n &=& A \Phi_{n-1} - \lambda_{n-1} \Phi_{n-1} \nonumber 
\end{eqnarray}
that is 
\[
\Phi_k = (A-\lambda_k I ) (A - \lambda_{k-1} I ) \cdots (A - \lambda_1 I)
\]
\end{lemma}

\begin{proof}
By induction. \\
Initialisation: 
\begin{eqnarray}
\Phi_1 &=& A \Phi_0 - \lambda_1 \Phi_0 = (A - \lambda_1 I) \nonumber \\
\Phi_2 &=& A \Phi_1 - \lambda_2 \Phi_1 = A (A - \lambda_1 I) - \lambda_2 I (A - \lambda_1 I) \nonumber \\
&=& (A - \lambda_2 I) (A - \lambda_1 I) \nonumber   
\end{eqnarray}
Induction:
Let suppose that $\Phi_k = (A - \lambda_k I)(A - \lambda_{k-1} I) \cdots (A  - \lambda_1 I)$.
Then the $k$th recursion step gives
\begin{eqnarray}
\Phi_{k+1} &=& A \Phi_k - \lambda_{k+1} \Phi_k \nonumber \\
&=& A \left[ (A - \lambda_k)(A - \lambda_{k-1}) \cdots (A - \lambda_1) \right] \nonumber \\
&& - \lambda_{k+1} \left[(A - \lambda_k I)(A - \lambda_{k-1} I) \cdots (A - \lambda_1 I)\right] \nonumber \\
&=& (A - \lambda_{k+1} I) \left[ (A - \lambda_k I)(A - \lambda_{k-1} I) \cdots (A  - \lambda_1 I) \right] 
\nonumber \\
&=& (A - \lambda_{k+1} I) (A - \lambda_{k} I) \ldots (A - \lambda_1 I) \nonumber 
\end{eqnarray}
which proves the induction step and the proof of the lemma.
\end{proof}

\begin{remark}
The procedure is simple and does not require computing directly $\Phi(A)$. The numerical behavior is roughly the same as with Ackermann's formula, hence, it is known to be numerically inacurate for ill-conditioned problems. Part of the explanation of this fact comes from the necessity to start with
$C$ and setting $K_0 = C$.
\end{remark}

One key feature in the upcoming modifications of Algorithm 1 is to develop the row vector $C$ as a product of matrices
\begin{equation}
\label{chain1}
C = C_{n-1} C_{n-2} \ldots C_1 
\end{equation}
where each matrix $C_{n-i}$ has $n-i$ columns and  $n-i+1$ rows, and the product up to and including $C_1$ cancels the spanning set of vectors of the controllability from $B$ to $A_{n-i} B$, that is, 
\[
C_{n-i+1} C_{n-i} \ldots C_1 \left( \begin{array}{cccc} B & AB & \ldots & A^{n-i} B  \end{array} \right) = 0 \qquad i=2,\ldots,n
\]

Although $C$ is unique, the chain of matrices $C_i$, $i=1,\ldots, n-1$ is not uniquely defined and finding the suitable chain is the important element for ensuring good numerical properties of the proposed method.
We will now illustrate the implication of the product chain of matrices $C_n$ $C_{n-1}$ $\ldots$ $C_n$ both (i) directly on the characteristic polynomial through suitable expansion and reformulation as an induction and (ii) applied to the above inductive factorization using the eigenvalues to be placed.  

\subsubsection{Chain of matrices applied to Ackermann's formula directly}
Let the desired characteristic polynomial be
\[
\Phi(s) = s^n + p_1 s^{n-1} + \ldots + p_{n-1} s + p_n
\]
which is the desired closed-loop polynomial after feedback $u = - K x$. 
Ackermann's formula becomes 
\begin{eqnarray}
K = C (A^n + p_1 A^{n-1} + \ldots + p_{n-1} A + p_n I) \label{facto1}
\end{eqnarray}

Let us for the moment introduce the expression (\ref{chain1}) into (\ref{facto1}) 

\begin{eqnarray}
K &=& C_{n-1} C_{n-2} \cdots C_2 C_1 A^n + C_{n-1} C_{n-2} \cdots C_2 C_1 A^{n-1} p_1 + \ldots \nonumber \\
&& + C_{n-1} C_{n-2} \cdots C_2 C_1 A p_{n-1} + C_{n-1} C_{n-2} \cdots C_2 C_1 p_n I \nonumber
\end{eqnarray}

\begin{eqnarray}
K &=& C_{n-1} C_{n-2} \cdots C_2 C_1 A^n + C_{n-1} ( C_{n-2} C_{n-3} \ldots C_2 C_1 A^{n-1} p_1  \nonumber \\
&& + C_{n-2} ( C_{n-3} C_{n-4} \cdots C_2 C_1 A^{n-2} p_2 + \ldots \nonumber \\
&& + C_2 ( A^2 p_{n-2} + C_1 (A p_{n-1} +  p_n I)) \ldots )) \label{Kexpr}
\end{eqnarray}
After defining
\begin{eqnarray}
A_1 &=& A \nonumber \\
A_k &=& C_{k-1} C_{k-2} \cdots C_2 C_1 A^k \qquad {k=2, \ldots, n}
\end{eqnarray}
the formula (\ref{Kexpr}) for the gain matrix becomes 
\begin{eqnarray}
K &=& A_n + C_{n-1} (A_{n-1} p_1 + C_{n-2} (A_{n-2} p_2 + \ldots + C_1 (A _1 p_{n-1} + p_n I) \ldots )) \nonumber \nonumber
\end{eqnarray}
which displays an inductive character starting with $K_0 = p_n I$ and 
$K_1 = A_1 p_{n-1} + C_1 K_0$ leading to the procedure
\begin{eqnarray}
K_0 &=& p_n I \nonumber \\
K_1 &=& C_1( A_1 p_{n-1} +  K_0) \nonumber \\
K_2 &=& C_2 (A_2 p_{n-2} +  K_1) \nonumber \\
\vdots && \vdots \nonumber \\
K_{n-1} &=& C_{n-1}(A_{n-1} p_1 +  K_{n-2}) \nonumber \\
K &=& A_n + K_{n-1} \nonumber 
\end{eqnarray}

\label{chainAckermannFactorization}

\subsubsection{Chain of matrices applied to the inductive factorization of Ackermann's formula}

Consider the inductive procedue (\ref{induction1}), and introduce the chain of matrices $C_{n-1} C_{n-2} \ldots C_1$ 
\begin{eqnarray}
K_0 &=& C_{n-1} C_{n-2} \ldots C_1  \nonumber \\
K_1 &=& K_0 A - \lambda_1 K_0 \nonumber \\
K_2 &=& K_1 A - \lambda_2 K_1 \nonumber \\
\vdots && \vdots \nonumber \\
K_{n-1} &=& K_{n-2} A - \lambda_{n-1} K_{n-2} \nonumber \\
K &=& K_{n-1} A - \lambda_n K_{n-1} \nonumber
\end{eqnarray}

The chain of matrices is distributed inside the induction steps without changing the result:
\begin{eqnarray}
K_0 &=& I  \nonumber \\
K_1 &=& K_0 A - \lambda_1 K_0 \nonumber \\
K_2 &=& C_1 K_1 A - \lambda_2 C_1 K_1 \nonumber \\
\vdots && \vdots \nonumber \\
K_{n-1} &=& C_{n-2} K_{n-2} A - \lambda_{n-1} C_{n-2} K_{n-2} \nonumber \\
K &=& C_{n-1} K_{n-1} A - \lambda_n C_{n-1} K_{n-1} \nonumber
\end{eqnarray}

This provides another formulation of the feedback computation, one that necessitates complex arithmetic when complex eigenvalues are required to be placed, and one that operates on the eigenvalues to be assigned.

 In case of complex conjugate eigenvalues, two steps of the above procedure step producing $K_k$ and $K_{k+1}$ can be grouped together (the eigenvalues are $\lambda_k = \lambda_{k+1}^{*} = \lambda \in \mathbb{C}$). 
\begin{eqnarray}
K_k &=& C_{k-1}K_{k-1}A - \lambda_k C_{k-2}K_{k-1}  \nonumber  \\
K_{k+1} &=& C_{k} K_k A - \lambda_{k+1} C_{k-1}K_k \nonumber
\end{eqnarray}
can be merged in one step using only real arithmetic
\begin{eqnarray}
K_{k+1} &=& C_{k-1} C_k K_k A^2 - (\lambda + \lambda^*) C_{k-1} C_k K_k A + \lambda \cdot \lambda^{*} C_{k-1} C_k K_k \nonumber
\end{eqnarray}
for which all three blocks $C_{k-1} C_k K_k A^2$, $C_{k-1} C_k K_k A$ and $C_{k-1} C_k K_k$ are real and computed from the block $C_{k-1} C_k K_k$.

\subsubsection{Computation of the chain of matrices $C_{n-1} C_{n-2} \ldots C_1$}

For the moment, we have assumed Ackermann's formula true, and the justification rests on the proof of validity of that formula. However, we do not have a particular way to compute the chain of matrices. The only thing required and sufficient for successful pole placement is that the product of the matrices appearing in the chain equals the last row of the controllability matrix, i.e.
\begin{equation}
C_{n-1} C_{n-2} \ldots C_{1} = e_n^T \, {\mathcal C}^{-1}
\end{equation}

We will show (in the Algorithms part) that the $C_i$ individually are proportional to anchors of specific algebroids, or in the linear algebra setting, specific projectors and commutators of matrices. Depending on the type of anchors (projectors), we will have different types of numerical properties--thanks to suitable scalings and orthogonal transformation-- and properties such as  whether we stay within the ring operations $(+,\cdot,\mathbb{Z})$ or not when the entries in both $A$ and $B$ belong to $\mathbb{Z}$.

\section{Intersection of affine hyperplanes}

Computing the single row vector $K$ that assigns the eigenvalues of $A - BK$ has a geometrical interpretation as the intersection of affine hyperplanes. The author has not found an explicit reference for this interpretation. Hence, the interpretation is considered as novel and will be given a complete exposition.

\subsection{The main intersection theorem}
Under some mild assumptions on the eigenvalues to be placed (distinct real eigenalues distinct from those of $A$), the following theorem will show that the gain vectors $K_i$ assigning one of these eigenvalues $\lambda_i$ (the other eigenvalues being arbitrary) constitute an affine hyperplane. Under the assumption of controllability, the orthogonals $n_i$ to the hyperplanes associated to each gain $K_i$ constitute a basis of $\mathbb{R}^n$ and the intersection of these affine hyperplanes intersect producing the single vector $K$ assigning the eigenvalues chosen.

\begin{theorem} 
\label{mainIntersectionTheorem} All the coefficients of the $B$ vector are supposed nonzero, i.e. $b_i \neq 0$, $i=1,\ldots n$. The eigenvalues to be placed $\lambda_i$ are supposed real and distinct from each other and not equal to the eigenvalues of the $A$ matrix. For each eigenvalue $\lambda_i$ compute $n$ gains given as
\begin{equation}
\label{pointHyperplane}
k_{ij} = 1/b_j (a_j^T - \lambda_i e_j^T) \qquad j=1,\ldots,n
\end{equation}
where $a_j^T$ is the $j$-th rowth of the $A$ matrix and $e_j$ is the $j$-th canonical basis.
The set of $k_{ij} \in \mathbb{R}^{1 \times n}$ for fixed $i$ and $j=1,\ldots,n$ constitute a $n\times n$ collection of row vectors. The following statements hold:
\begin{enumerate}
\item[i)] Formula (\ref{pointHyperplane}) provides vectors $k_{ij}$ $j=1,\ldots,n$ that assigns the eigenvalue $\lambda_i$, that is,
\[
| \lambda I - A + B k_{ij}|= p_j(\lambda) (\lambda - \lambda_i)
\] 
where $p_j(\lambda)$ is some polynomial of degree $n-1$.
\item[ii)]
The  $k_{ij}$, $j=1,\ldots,n$ are linearly independent. This holds for each $i = 1, \ldots, n$. Fixing $i$ and $j$, the differences $k_{ij} - k_{ij^{'}}$, $j^{'}=1, \ldots, n$, $j^{'} \neq j$ consitute a basis of a $n-1$ vector space which is a basis of the hyperplane associated with $\lambda_i$ parallel to the affine hyperplane associated with $\lambda_i$. The affine hyperplane associated with $\lambda_i$ is the one parallel that contains the points  $k_{ij}$, $j=1, \ldots, n$.
\item[iii)]
Any vector connecting the origin to a point on the affine hyperplane associated with $\lambda_i$  assigns the eigenvalue $\lambda_i$.
\item[iv)]
The hyperplanes intersect at a point giving the vector $K$, if, and only if, the system is controllable. The point of intersection is the end point of the vector $K$ such that
\[
| \lambda I - A + B K|= \prod_{i=1}^n (\lambda - \lambda_i)
\] 
\end{enumerate}
\end{theorem}

Proof:  To show i), we can consider without loss of generality a given $i$ and $j$. Expressing the matrices appearing according to their rows, and expanding $B k_{ij}$ using the definition (\ref{pointHyperplane}) gives
\begin{eqnarray}
&& \lambda_i I - A + B k_{ij}  = \nonumber \\
&&  \left( \begin{array}{c}
\lambda_i e_1^T \\
\vdots \\  \lambda_i e_{j-1}^T \\ \lambda_i e_j^T \\ \lambda_i  e_{j+1}^T \\ \vdots \\
\lambda_i e_n^T \end{array} \right)
 - \left( \begin{array}{c}
a_1^T \\  \vdots \\ a_{j-1}^T \\ a_j^T \\ a_{j+1}^T \\ \vdots \\  a_n^T \end{array} \right) + \left( \begin{array}{c}
\frac{b_1}{b_j}  (a_{j}^T- \lambda_i e_j^T) \\
\vdots \\
\frac{b_{j-1}}{b_j} (a_{j}^T- \lambda_i e_j^T) \\
a_j^T - \lambda_i e_j^T \\
\frac{b_{j+1}}{b_j} (a_{j}^T- \lambda_i e_j^T)  \\
\vdots \\
\frac{b_n}{b_j} (a_{j}^T- \lambda_i e_j^T)  \end{array} \right)
= \left( \begin{array}{c}
** \\
\vdots \\
** \\
0 \\
** \\
\vdots \\
** \end{array} \right) \nonumber
\end{eqnarray}
where $**$ stands for an arbitrary row of not necessarily zero values. Because of a row of zeros appearing on the $j$-th row means that taking the determinant gives
\[
 | \lambda_i I - A + B k_{ij} | = 0
\]
which proves that the eigenvalue $\lambda_i$ satisfies the charecteristic polynomial of $A - B k_{ij}$ and  is therefore indeed assigned, thus proving the statement i). 

The independence of the $k_{ij}$ appearing in  ii) can be proved by looking at their expression given by (\ref{pointHyperplane}). By rescaling each vector by a nonzero scalar does not change their linear dependency. Hence, fixing $i$, and considering $b_j k_{ij}$ instead of $k_{ij}$ gives
\[
\left| \begin{array}{c} b_1 \, k_{i1} \\ b_2\, k_{i2} \\ \vdots \\ b_n\, k_{in}
\end{array} \right| = \left| \begin{array}{c} a_1^T - \lambda_i\, e_1^T \\
a_2^T - \lambda_i \, e_2^T \\ \vdots \\
a_n^T - \lambda_i \, e_n^T \end{array} \right| = | A - \lambda_i I| \neq 0
\]
The determinant is nonzero because $\lambda_i$ is not an eigenvalue of $A$. This is true for all the $i = 1,\ldots,n$.
Now, since the $k_{ij}, j=1,\ldots,n$ are linearly independent, the 
vectors $k_{ij} - k_{ij^{'}}$ for $j \neq j^{'}$, $j,j^{'} = 1,\ldots, n$ are linearly independent and we can define the affine Hyperplane $i$ to be the affine hyperplane containing the end point of $k_{ij}$ and spanned by the $n-1$ vectors $k_{ij} - k_{ij^{'}}$ for $j \neq j^{'}$, $j \neq i$, $j^{'} \neq i$, and $j,j^{'} = 1,\ldots, n$.

To prove iii), let us start by showing that $k_{11} + \alpha(k_{11} - k_{12})$ sets one of the eigenvalues of $A - B (k_{11} + \alpha(k_{11}- k_{12}))$ to $\lambda_1$, for all choices of $\alpha$.
Setting 
\[
v^T = \frac{b_2}{b_1} \left(a_1^T - \lambda_1 e_1^T\right) - a_2^T + \lambda_1 e_2^T
\]
Formula (\ref{pointHyperplane}) gives
\begin{eqnarray}
&&\lambda_1 I - A + B k_{11} + \alpha B(k_{11}- k_{12}) = \nonumber \\
&&\left( \begin{array}{c}
\lambda_1 e_1^T - a_1^T + a_1^T - \lambda_1 e_1^T + \alpha \frac{b_1}{b_2} v^T \\
\lambda_1 e_2^T - a_2^T + \frac{b_2}{b_1}(a_1^T - \lambda_1 e_1^T) + \alpha v^T \\
\lambda_1 e_3^T - a_3^T + \frac{b_3}{b_1}(a_1^T - \lambda_1 e_1^T) + \alpha \frac{b_3}{b_2} v^T \\
\vdots \\
\lambda_1 e_n^T - a_n^T + \frac{b_n}{b_1}(a_1^T - \lambda_1 e_1^T) + \alpha \frac{b_n}{b_2} v^T
\end{array} \right) \label{matrixSecond}
\end{eqnarray}
Taking the determinant of this matrix gives zero, because the second row of (\ref{matrixSecond}) is $v^T + \alpha v^T$ and and the first row
is $0 + \frac{b_1}{b_2} \alpha v^T$ and all sucessive rows are obtained from the rows of $A - \lambda_1 I - B k_{11}$
by adding a suitable multiple of $v^T$.
Indeed, if $v^T \in \text{span}\{ v_1^T, \ldots, v_{n-1}^T\}$ then
\[
\left| \left( \begin{array}{c} 0 \\ v_1^T \\ \vdots \\ v_n^T \end{array} \right) +
\left(\begin{array}{c} \beta_1 v^T \\ \beta_2 v^T \\ \vdots \\ \beta_n v^T \end{array} \right) \right| = 0
\] no matter the values of $\beta_1, \ldots \beta_n \in \mathbb{R}$.
Now, setting $k_{11}^{'} = k_{11} + \alpha(k_{11} - k_{12})$ the process can be repeated to show
\[
| \lambda_1 I - A - B k_{11}^{'} + \alpha^{'}( k_{11}^{'} - k_{13}^{'}) |= 0
\]
and proceeding inductively we finally arrive at 
\[
\left| \lambda_1 I - A - B k_{11}^{'} + \sum_{i=2}^n \alpha_i (k_{11}- k_{1i}) \right|= 0
\] which proves statement iii) once the whole process is repeated again for each eigenvalue to be placed (i.e. $k_{2j}$ for $\lambda_2$, $k_{3j}$ for $\lambda_3$, $\ldots$, 
$\lambda_n$ for $\lambda_{nj}$, $j=1, \ldots, n$ in every case).

To prove iv) we start by showing that if the system is controllable then the affine hyperplanes intersect to give $K$. By assumption and using Ackermann's formula \cite{Ackermann}, we can assign any distinct real eigenvalues $\lambda_i, i=1, \ldots, n$, distinct from those of $A$. We can use formula (\ref{pointHyperplane}) to find the $n$ affine hyperplanes, each one assigning one of the eigenvalues $\lambda_i$ chosen. Since each affine hyperplane has a supporting vector space of dimension $n-1$, each affine hyperplane, say $H_i$ associated with $\lambda_i$, must pass through $K$ because $K$ assigns $\lambda_i$. This means that $K$ is common to every affine hyperplane $H_i$, $i=1,\ldots n$. The $n$ affine hyperplanes intersect giving $K$.

To prove the reciprocal, we use the contraposition principle. We will show that if the system is not controllable then the affine hyperplanes cannot intersect. This will be shown by displaying a contradiction. If the system is not controllable, the Popov-Belevitch-Hautus test states that  
the matrix \[ [A - \tilde \lambda I \qquad B ]\]
loses rank for at least one eigenvalue $\tilde \lambda$ of $A$, see \cite{Chen}, \cite{Kailath}. This eigenvalue remains an eigenvalue of
\[A - B K\]
no matter the choice of $K$. (Indeed, decompose the system using similarity transform $P$ into the controllable and uncontrollable spaces. The transformed $PB$ is zero for the uncontrollable states. The uncontrollable states are governed by a submatrix $A_{\text{nc}}$ that has $\tilde \lambda$ as an eigenvalue. Hence $\tilde \lambda$ cannot be changed. Details are in \cite{Chen} or \cite{Kailath}, for example.)

Choose $n$ real and distinct eigenvalues $\lambda_i$ distinct from those of $A$ (and in particular distinct from $\tilde \lambda$). Using point ii) of the Theorem, it is possible to construct $n$ affine hyperplanes, each one assigning one of the $\lambda_i$, $i = 1, \ldots n$. If these affine hyperplanes intersect, the intersection gives a gain vector $\bar K$ that assigns eigenvalues different from $\tilde\lambda$ which is a contradiction. Therefore, the affine hyperplanes do not intersect.

To show iv), we use the assumption of controllability, and using iii), it is guaranteed that the hyperplanes constructed in ii) intersect giving the vector $K$ that assigns arbitrary real and distinct eigenvalues, and distinct from those of $A$.

This concludes the proof of the theorem.

\subsection{Examples illustrating the main intersection theorem}

Two examples will illustrate the theorem and points appearing in the proof of Theorem \ref{mainIntersectionTheorem}.
The first example is a controllable one. The second is an uncontrollable one.

The first example is
\[
A = \left( \begin{array}{ccc} 1 &3 &5 \\ 7 &13 &17 \\ 1 &1 &1 \end{array} \right) \qquad B = \left( \begin{array}{c} 1 \\ 1 \\ 1 \end{array} \right)
\]
It is controllable since
\[
| B \, AB \,A^2 B | = 352 \neq 0
\]
Without loss of generality, the eigenvalues chosen are $\lambda_1 = -1$, $\lambda_2 = -2$, $\lambda_3 = -3$. The script for this example can be found in Section \ref{exampleIntersection}.
For $\lambda_1=-1$ we get the following row vectors of gains using formula (\ref{pointHyperplane})
\begin{eqnarray}
k_{11} &=& \left(  \begin{array}{ccc} 2 &  3 & 5 \end{array} \right) \nonumber \\
k_{12} &=& \left( \begin{array}{ccc} 7 & 14 & 17 \end{array} \right) \nonumber \\
k_{13} &=& \left( \begin{array}{ccc} 1 & 1 & 2 \end{array} \right) \nonumber
\end{eqnarray}
They are linearly independent since
\begin{equation}
\label{detEx}
\left| \begin{array}{ccc} 2 & 3 & 5 \\ 7 & 14 & 17 \\ 1 & 1 & 2 \end{array}
\right| = - 4 
\end{equation}
We check the predictions of the theorem, that each $k_{1i}$, $i=1,2,3$, assigns one eigenvalue at $\lambda_1 = -1$:
\begin{eqnarray}
| \lambda I - A + B k_{11} | &=& (\lambda +1) (\lambda - 8) (\lambda+2) \nonumber \\
| \lambda I - A + B k_{12} | &=& (\lambda + 1) (\lambda^2 + 22 \lambda + 24) \nonumber \\
| \lambda I - A + B k_{13} | &=& (\lambda + 1) (\lambda^2 - 12 \lambda - 12) \nonumber 
\end{eqnarray}
Similar results hold for $k_{2j}$ with eigenvalue $\lambda_2 = -2$ and $k_{3j}$ with eigenvalue $\lambda_3 = -3$.
The affine hyperplane 1 has the following equation (the details and a reference for determinental intersections are in Section \ref{exampleIntersection}).

\[
0 = -4 + 9 \, x - 3 \,y - z 
\]
The affine hyperpane 2 has the following equation
\[
0 = 32 - 0 \, x + 16 \, y - 16 \, z
\]
The affine hyperplane 3 has the following equation
\[
0 = 110 - 11 \, x + 33 \, y - 33 \, z
\]
The three affine hyperplanes intersect since
\[
\left| \begin{array}{ccc} 9 &-3 &-1 \\
0 & 16 & -16 \\
-11 & 33 & -33 \end{array} \right| = - 704
\]
The point of intersection is
\begin{eqnarray}
K &=& \frac{1}{-704} \left( \begin{array}{ccc} -2816 & -5280 & -6688 \end{array} \right) \nonumber \\
&=& \left( \begin{array}{ccc} 4 & \frac{15}{2} & \frac{19}{2} \end{array} \right) \nonumber
\end{eqnarray}
And a direct computation gives 
\[
| \lambda I - A + B K | = (\lambda+1)(\lambda+2)(\lambda+3)
\]
confirming the correct placement of the eigenvalues.
Notice that due to the choice of $1$ in the $B$ matrix all operations except the last division could be performed in the ring $(+, \cdot , \mathbb{Z})$ since all entries of $A$ were in $\mathbb{Z}$ and those of $B$ were units, and the intersection and computation of the intersection relied on determinants.

The second example is 
\[
A = \left( \begin{array}{ccc} 6 & 4 & -9 \\
5 & 2 & -6 \\
0 & 0 & 1 \end{array} \right ) \qquad B = \left( \begin{array}{c} 1\\ 1 \\ 1 \end{array} \right)
\]
It is an uncontrollable system since
\[
| B \, AB \,A^2 B | = 0
\]
There is no difficulty in assigning individually each $\lambda_1 = -1$, $\lambda_2 = -2$ and $\lambda_3 = -3$ and the hyperplanes are well defined. However, they do not intersect. Indeed the equations of the hyperplanes can be computed as in the previous example (see Section \ref{exampleIntersection})
\begin{eqnarray}
0 &=& 2 - x - y -z \nonumber \\
0 &=& 36 -12\, x -12\, y -12 \, z \nonumber \\
0 &=& 100 -25 \,x -25 \, y - 25 \, z
\end{eqnarray}
They do not intersect since
\[
\left| \begin{array}{ccc} -1 & -1 & -1 \\  -12 & -12 & -12 \\ -25 & -25 & -25 \end{array} \right| = 0 
\]
The affine hyperplanes are all parallel to each other.

\label{examples3}

\part{Algorithms}

\section{Algorithm based on the intersection of affine hyperplanes}

We will illustrate and/or prove three ways of determining the intersection of the affine hyperplanes

\begin{itemize}
\item Finding the intersection using determinants. 
\item Successive projections by sliding along the hyperplanes.
\item Setting one eigenvalue at a time and quotienting out the result obtained by quotienting into the hyperplane associated with the eigenvalue placed. We will call this method the affine hyperplane intersections using quotients, or the first algebroid method.
\end{itemize}

\subsection{Intersection of affine hyperplanes using determinants}

This method  has been illustrated in the previous section and we will not go further in this direction because we will provide a much more efficient way of computing the gain using only ring arithmetic later on.
It is quite straightforward to generalize the example to $n$ dimensions (see the script for the 3 dimensional example).

\subsection{Intersection of affine hyperplanes  using successive projections}

A technique to find the intersection of the $n$ affine hyperplanes
is to start by computing orthogonal vectors $n_i$  to each hyperplane $H_i$, $i=1,\ldots,n$. Then starting from an arbitrary
$K_1$ that assigns $\lambda_1$, the procedure is to slide perpendicularly to $n_1$ until $H_2$ is reached. Then repeat
the operation by sliding perpendicularly to both $n_1$ and $n_2$ until $H_3$ is reached. In this way, it is guaranteed that $\lambda_1$ and $\lambda_2$ remain assigned while sliding. Proceeding inductively, the step $j$ is reached, and then slide perpendicularly to the vectors $n_1, n_2, \dots, n_j$ until the affine hyperplane $H_{j+1}$ is reached. The algorithm stops when reaching the affine hyperplane $H_n$. The final destination point gives the intersection and the feedback gain $K$.

Before giving the method, we give a lemma to find the orthogonal vectors $n_i$.

\begin{lemma}
The orthogonal $n_i$ is the solution to 
\begin{equation}
\label{normals}
(A - \lambda_i I) n_i = B
\end{equation}
\end{lemma}

Proof:
Examining (\ref{pointHyperplane}) gives
\begin{eqnarray}
(A - \lambda_i I) n_i  = \left( \begin{array}{c}
b_1 k_{i1} \\
\vdots \\
b_n k_{in}  \end{array} \right) n_i  \label{contin}
\end{eqnarray}
Now setting
\[
n_i = \left( \begin{array}{c}
k_{i1} \\
\vdots \\
k_{in} 
\end{array} \right)^{-1} \left( \begin{array}{c} 1 \\ \vdots \\ 1 \end{array} \right)
\]
which is well defined because by assumption the $k_{ij}$ are linearly independent for a each fixed $i$, 
gives the solution mentioned in the lemma, that is,
\[
(A - \lambda_i I) n_i = B
\]
Notice first that $k_{ij} n_i = 1$ for all $j = 1, \ldots, n$ by construction.
Indeed, 
\[
\left( \begin{array}{c}
k_{i1} \\
\vdots \\
k_{1n} 
\end{array} \right) n_i = \left( \begin{array}{c} 
k_{i1} \\
\vdots \\
k_{in} \end{array} \right)
\left( \begin{array}{c}
k_{i1} \\
\vdots \\
k_{in} \end{array} \right)^{-1} \left( \begin{array}{c} 1 \\ \vdots \\ 1
\end{array} \right) = \left( \begin{array}{c} 1 \\ \vdots \\ 1 \end{array} \right)
\]
Finally, multiplying $k_{ij} n_i = 1$ on both sides by $b_j$ gives
\[
b_j k_{ij} n_i = b_j
\] and shows that (\ref{contin}) equals $B$. Since $k_{ij} n_i = 1$, one has 
$(k_{ip} - k_{ir}) n_i = 0$, no matter $p$, and $q$ in $1, \ldots n$ and therefore $n_i$ is orthogonal to the hyperplane $i$ that $k_{ij}$ define for a fixed $i$. This completes the proof of the lemma.

The procedure described earlier (successive slidings and intersection) is now given. Finding the intersections is preformed by suitable projections. The $G_i$ are the projectors, and they have an interpretation in the theoretical Section \ref{anchors}, as part of a commutator and anchor in the Algebroid setting. 

\begin{enumerate}
\item Compute the orthogonal vectors to the affine hyperplanes $n_1$, $n_2$, $\ldots$, $n_n$ according to the relations (\ref{normals}). The vectors $n_i$ --- one for each hyperplane, $i=1,\ldots,n$  --- are obviously not normalized and are given as column vectors of size $n \times 1$.
\item  Compute for each hyperplane one point belonging to it --- according to  for example formula  (\ref{pointHyperplane}) --- so as to obtain row vectors  $K_1$, $K_2$, $K_3$, $\ldots$, $K_n$ with the index $i$ designating the  hyperplane where the row vector $K_i$ of size $1 \times n$ belongs. The $K_i$'s are used in Step 5, below.
\item
Set $n_{i0} := n_i$ for $i = 1, \ldots, n$ to start the following iteration steps:
\begin{eqnarray}
n_{i1} &:=& \left(I - \frac{n_{10} n_{10}^T}{ n_{10}^T n_{10}} \right) n_{i0} \qquad i = 2,\ldots,n \nonumber  \\
n_{i2} &:=& \left(I - \frac{n_{21} n_{20}^T}{n_{20}^T n_{21}} \right) n_{i1} \qquad i=3, \ldots,n \nonumber \\
\vdots && \vdots \nonumber  \\
n_{ik} &:=& \left(I  - \frac{n_{k,k-1} n_{k0}^T}{n_{k0}^T n_{k,k-1}} \right) n_{i,k-1} \qquad i=k+1, \ldots, n \nonumber \\
\vdots && \vdots \nonumber \\
n_{n,n-1} &:=& \left( I - \frac{ n_{n-1,n-2} n_{n-1,0}^T}{n_{n-1,0}^T n_{n-1,n-2}} \right) n_{n,n-2} \nonumber 
\end{eqnarray}
\item  Rename the rank-one matrices obtained in Step 3 according to:
\begin{eqnarray}
G_1 &:=& \frac{n_{21} n_{20}^T}{n_{20}^T n_{21}} \nonumber \\
G_2 &:=& \frac{n_{32} n_{30}^T}{n_{30}^T n_{32}} \nonumber \\
\vdots && \vdots \nonumber \\
G_{n-1} &:=& \frac{n_{n,n-1} n_{n,0}^T}{n_{n,0}^T n_{n,n-1}} \nonumber  
\end{eqnarray}
\item 
\begin{eqnarray}
\gamma_1 &:=& K_1 \nonumber \\
\gamma_2 &:=& \gamma_1 - (\gamma_1 - K_2) G_1^T \nonumber \\
\gamma_3 &:=& \gamma_2 - (\gamma_2 - K_3) G_2^T \nonumber \\
\vdots && \vdots \nonumber \\
\gamma_n &=& \gamma_{n-1} - (\gamma_{n-1} - K_n) G_{n-1}^T \nonumber 
\end{eqnarray}
\item \[ K := \gamma_n \]
\end{enumerate}

To understand the above alorithm, one should prove that the $\gamma_i$, $i = 1, \ldots, n$ are the result of the sliding
and intersection mechanism with the affine Hyperplanes described earlier, that is $\gamma_i$ assigns the eigenvalues
$\lambda_1, \ldots, \lambda_i$. We will show this property by recurrence.

We first proceed with the initalisation step of the recurrence.
Recall that the hyperplane of gains for the first eigenvalue $\lambda_1$ is the set  of points corresponding to all possible $K_1$'s that assign the eigenvalue $\lambda_1$. The row  vector $K_1$ is the vector from the origin  to one of the points of the first affine hyperplane.  Then in a similar fashion, $K_2$ is one of the possible vectors assigning the second real eigenvalue $\lambda_2$. All such $K_2$ vectors connect the origin to one of the points of the second affine hyperplane. These affine hyperplanes will be designated explicitely with a capital 'H', affine Hyperplane 1 and affine Hyperplane 2. We will also drop the adjective 'affine' to simplify notation, and call them  Hyperplane 1 and Hyperplane 2.

Therefore, starting from $K_1= \gamma_1$, the row vector $\gamma_2$ connects the origin to one of the points lying on the intersection of Hyperplanes 1 \& 2. We state this fact as a lemma.

\begin{lemma}
$\gamma_2$ assigns the real eigenvalues $\lambda_1$ and $\lambda_2$. The row vector $\gamma_2$ connects the origin to a point lying on the intersection of the affine hyperplane associated to the gains $K_1$ (Hyperplane 1) and the affine hyperplane of associated to the gains $K_2$ (Hyperplane 2).
\end{lemma} 

\paragraph{Proof: } It is sufficient to prove the following two assertions 
\begin{enumerate}
\item $K_1-\gamma_2$ is in a hyperplane parallel to Hyperplane 1.
\item $K_2 -\gamma_2$ is in a hyperplane parallel to Hyperplane 2.
\end{enumerate}
Once those two assertions are established the conclusion of the lemma follows directly, since the first assertion implies that the end point of $\gamma_2$ is in Hyperplane 1 (the hyperplane among all parallel hyperplanes sharing the same normal  must contain the end point of $K_1$, since the end point of $\gamma_2$ is the prolongation of the end point of $K_1$ by a vector $\gamma_2 - K_1$ orthogonal to all hyperplanes parallel to  Hyperplane 1). The second assertion implies that among all parallel hyperplanes sharing the common normal, the end point of $\gamma_2$ lies on the same hyperplane as the end point of $K_2$, that is, it lies on Hyperplane 2. The argument is similar as for the first assertion ($\gamma_2$ is constructed as the prolongation of the end point of $K_2$--- lying on Hyperplane 2 ---  by a vector orthogonal to normal of Hyperplane 2).

Hence let us start with the first assertion recalling the definition of $G_1= \frac{n_{21}n_{20}^T}{n_{20}^T n_{21}}$. The assertion follows from the fact that the 1st orthogonal is $n_{10}$ 
\begin{eqnarray}
G_1^T n_{10} &=& \frac{n_{20} n_{21}^T}{n_{21}^T n_{20}} n_{10} = \frac{n_{20}}{n_{21}^T n_{20}} \, n_{21}^T n_{10} \nonumber \\
n_{21} &=& \left(I - \frac{n_{10} n_{10}^T}{n_{10}^T n_{10}} \right) n_{10} \nonumber \\
n_{21}^T &=& n_{10}^T \left(I - \frac{n_{10} n_{10}^T}{n_{10}^T n_{10}} \right)  \nonumber \\
G_1^T n_{10} &=& \frac{n_{20}}{n_{21}^T n_{20}} n_{10}^T \left( I -
\frac{n_{10} n_{10}^T}{n_{10}^T n_{10}} \right) n_{10} \nonumber \\
&=& \frac{n_{20}}{n_{21}^T n_{20}} \left( n_{10}^T n_{10} - \frac{ n_{10}^T n_{10} n_{10}^T n_{10}}{n_{10}^T n_{10}} \right) \nonumber \\
&=& 0 \nonumber \\
(\gamma_2 - K_1) n_{10} &=& - (\gamma_1 - K_2) G_1^T n_{10}
= 0 \nonumber 
\end{eqnarray}
The second assertion follows similarly using the second orthogonal $n_{20}$
\begin{eqnarray}
(K_2 - \gamma_2) n_{20} &=& (K_2 - \gamma_1 + (\gamma_1 - K_2) G_1^T) n_{20} \nonumber \\
&=& (K_2 - \gamma_1) \left( I - G_1^T \right) n_{20} \nonumber \\
&=& (K_2 - \gamma_1) \left( I - \frac{n_{20} n_{21}^T}{n_{21}^T n_{20}} \right) n_{20} \nonumber \\
&=& (K_2 - \gamma_2) \left( n_{20} - n_{20} \frac{n_{21}^T  n_{20}}{n_{21}^T n_{20}} \right) \nonumber  \\
&=& 0  \nonumber
\end{eqnarray}
This concludes the proof of the initialisation stage of the recurrence.

The inductive step follows a similar reasoning and it is left to the reader.

The script for the 3 dimensional controllable example given in Section \ref{examples3} is given in Section \ref{exampleAlgebroid}.

{\bf Remark: } There are many arbitrary choices to be made, such as the values of the $K_i$. A further study is required to understand the impact of these choices on the final result due to error arithmetic. 

\subsection{The first algebroid method, quotienting into the hyperplanes}

The previous intersection method, although correct from a purely mathematical point of view, suffers from a few defects concerning the implementation, and in particular the numerical ones, mainly for the following reasons.
\begin{itemize}
\item It is not clear how to select the $K_i$, $i = 1, \ldots n$.
\item The sliding mechanism has to be precise in the computations of orthogonal directions to make sure that the previous eigenvalues placed remain at their position when sliding along the new direction. 
\item It is not obvious numerically which affine hyperplane to select when sliding.
\end{itemize}

The method proposed in this section circumvents the first issue by providing a well-defined point on the affine hyperplane. The idea is to use the orientation of the $B$ vector to find the unique point on the affine hyperplane. In this way, we also generalize the application to systems for which some entries of the $B$ vector are zeros.

It circumvents the second issue by performing a quotient operation using the  orthogonal anchor associated with the normal of the current hyperplane. This reduces the dimension one by one. The accuracy of this method will be shown through critical examples in the numerical experiments section.

It leaves the last issue open, that is, in which order to select the quotienting directions (i.e. which of the affine hyperplanes to select). This issue will be illustrated through examples in the numerical experiments  section.

Hence the method relies on reducing one dimension at a time, and placing one eigenvalue one at a time, at a typical stage of the descending part of the algorithm. 
The gain vector is constructed in the second part of the algorithm by proceeding backwards along ascending stages increasing the dimension one by one.

A typical descending stage in the first part of the algorithm  is decomposed in the following steps:
\begin{enumerate}
\item Compute the orthogonal direction to the current $B$ vector. Perform a qr decomposition of B
to obtain the anchor $Q^T$
\[
\left( \begin{array}{c} q^T \\ Q^T \end{array} \right) \qquad Q^T B = 0 \qquad q^T B \neq 0
\]
\item Using these directions, compute an orthogonal basis $Q_s^T$ of the current hyperplane (which will be the anchor for quotienting into the hyperplane). To achieve this,
let $\lambda_i$ denote the current eigevalue to be placed compute a qr decomposition of 
\[(A - \lambda_i \, I)^T \, Q\]
denoted \[ \left( \begin{array}{c}  Q_s^T \\ q_s^T \end{array} \right) \]
where $q_s^T$ is associated with the null values of the qr decomposition (i.e. the orthonormal to the hyperplane).
\item Place the eigenvalue corresponding the the hyperplane using the $B$ vector orientation. This determines a point on the hyperplane which is a $k$ vector.
\begin{equation}
\label{Bhat}
k = \frac{{\hat B}^T}{\| B \|} ( A - \lambda_i I)
\end{equation}
\item Project this $k$ vector so that the new $k$ vector is situated on the orthogonal line to the current affine hyperplane leading to the origin.
\[
k_o = k \, q_s q_s^T
\]
\item Perform a quotient operation using the generators of the hyperplane to define the anchor to obtain the next $A$ matrix and the next $B$ vector. 
\[ \bar A = Q_s^T A Q_s \qquad \bar B = Q_s^T B \]
\end{enumerate}

Repeat all these steps until there is no more dimensions.

\label{algebroidDescription}
The computer implementation for the general $n$ dimensional case is furnished in Section \ref{algebroidAlgorithm}. The steps corresponding to the first descending part of the algorithms are reapeated here for convenience, with some comments added.
\begin{verbatim}
n = length(B);
Ab = A;
Bb = B;
ani = zeros(n*(n-1),n);
kos = zeros(n,n);

for i=1:n-1

% compute the orthogonal direction to B (Q^T and q^T):
[qa,ra]=qr(Bb);

% compute the generators of the hyperplane:
[qs,rs]=qr((qa(2:end,:)*(Ab-VP(i)*eye(n-i+1)))'); 

% get the unit normal vector to the hyperplane:
koh = qs'(end,:); 

% perform Steps 3 and 4 in one step:
ko = Bb'/(Bb'*Bb)*(Ab - VP(i)*eye(n-i+1))*koh'*koh; 

% the anchor is Q_s^T of the theory (gen. of hyperplane)
anb = qs'(1:end-1,:);

% storage for later use (for the ascending part of algorithm):
ani((i-1)*n+1: (i-1)*n+n-i, 1: n-i+1) = anb;
kos(i,1:n-i+1)=ko; 

% quotienting A and B into the hyperplane:
Ab = anb*Ab*anb'; 
Bb = anb*Bb;

end;
\end{verbatim}

Notice in the code that $\texttt{koh}$ is  $q_s^T$ of the description of the descending steps.

Let us check what the code produces with the 3-dimensional example (our "fil rouge"):
\[ A = \left( \begin{array}{ccc} 1 & 3 & 5 \\ 7 & 13 & 17 \\ 1 & 1 & 1 \end{array} \right) \quad 
B = \left( \begin{array}{c} 1 \\ 1 \\ 1 \end{array} \right) \]
Running the above script on this examples produces using SysQuake 6.5
with $-1$, $-2$ and $-3$ as the  eigenvalues to be placed, in that order
\begin{verbatim}
koh =
   -0.9435    0.3145    0.1048
ko =
    0.3956   -0.1319   -0.0440
anb =
   -0.2923   -0.6401   -0.7105
   -0.1563   -0.7010    0.6959
Ab =
   17.2209   -1.2808
   15.4917   -1.4406
Bb =
   -1.6429
   -0.1614
koh =
   -0.0503    0.9987
ko =
   -0.0688    1.3654
anb =
   -0.9987   -0.0503
Ab =
   17.8876
Bb =
    1.6490
\end{verbatim}
The first \texttt{ko} places the eigenvalue at $-1$ for the orignal $A$ and $B$.
\begin{verbatim}
> ko1 = kos(1,1:n)
ko1 =
    0.3956   -0.1319   -0.0440
> eig(A-B*ko1)
ans =
   16.0889
   -1.0000
   -0.3087
\end{verbatim}
the other two eigenvalues are still not well defined at this stage.
The second \texttt{ko} places the eigenvalue $-2$ for the quotient system $\bar A$ and $\bar B$
The first quotient system is
\begin{verbatim}
Ab =
   17.2209   -1.2808
   15.4917   -1.4406
Bb =
   -1.6429
   -0.1614
\end{verbatim}
and one can check that with the second \texttt{ko}
\begin{verbatim}
> ko
ko =
   -0.0688    1.3654
> eig(Ab-Bb*ko)
ans =
   17.8876
   -2.0000
\end{verbatim}
showing that $-2$ has been set and the remaining eigenvalue still to be set.

\vskip 0.2cm
It's now time to provide the bottom-up part of the algorithm. The idea is to merge the supposedly well-defined gain vector on the quotient system with the system before taking the quotient, where only one eigenvalue is defined correctly, the one associated with the current affine hyperplane, the others still being undefined. The trick is that the gain of the quotient system is in one-to-one correspondence with a point on the hyperplane before taking the quotient. This point defines the same eigenvalues as those of the quotient system. This means that all eigenvalues are set properly if those of the quotient are set properly.

\begin{enumerate}
\item Set the scalar value from the last stage of the descending stage to 
\[K_n = \frac{1}{\bar B} (\bar A - \lambda_n)\]
\item Iterate for $i = n-1$ to $i =1$
\begin{equation}
\label{pull}
 K_{i} = k_{o,i} + K_{i+1} \, Q_{s,i}^T 
\end{equation}
\end{enumerate}
or in code
\begin{verbatim}
K = (Ab-VP(n))/Bb;
for i=n-1:-1:1
anb = ani((i-1)*n+1:(i-1)*n+n-i,1:n-i+1);
ko = kos(i,1:n-i+1);
K = ko + K*anb;
end;
\end{verbatim}
Running the above code on our example gives
\begin{verbatim}
Ab =
   17.8876
Bb =
    1.6490
K =
   12.6671
anb =
   -0.9987   -0.0503
ko =
   -0.0688    1.3654
K =
  -12.7198    0.7282
anb =
   -0.2923   -0.6401   -0.7105
   -0.1563   -0.7010    0.6959
ko =
    0.3956   -0.1319   -0.0440
K =
    4.0000    7.5000    9.5000
\end{verbatim}
which gives the correct gain $K$ such that 
\[| A - B K - \lambda I | = (\lambda +1)(\lambda + 2)(\lambda + 3) \].

To be complete the correctness of the algorithm,  we will need to fill the missing justifications:
\begin{itemize}
\item The formula for the gain vector (\ref{Bhat}) assigns the eigenvalue $\lambda_i$.
\item $Q_s^T$ are orthogonal generators of the affine hyperplane associated with $\lambda_i$.\item The eigenvalues of the quotient are carried over to the the system prior to taking the quotient when using formula (\ref{pull}). This formula leaves unchanged the current well placed eigenvalue $\lambda_i$.
\item Show that the algorithm works for multiple identical real eigenvalues.
\item Explain what happens when there is a loss of controllability. 
\end{itemize}
\subsubsection{The formula for the gain vector (\ref{Bhat}) assigns the
eigenvalue $\lambda_i$}

The proof follows a similar direct computation as for the proof of formula (\ref{pointHyperplane}) using a direct expansion. Setting
\[
k = \frac{{\hat B}^T}{\| B\|} (A - \lambda_i \, I)
\]
and then 
\begin{eqnarray}
| \lambda_i I - A + B k | &=& | \lambda_i I -  A + {\hat B} {{\hat B}^T} (A - \lambda_i \, I) | \nonumber \\
&=& | (I - {\hat B} {\hat B}^T) (-A + \lambda_i I) | \nonumber \\
&=& | I - {\hat B} {\hat B}^T | | -A + \lambda_i I | \nonumber
\end{eqnarray}
The second factor is non zero because $\lambda_i$ is by hypothesis not an eigevalue of $A$. The first factor is zero and this is easily shown by choosing a complementary set of $n-1$ vectors $\hat B_i$, $i=1, \ldots, n-1$ to constitute an orthogonal basis so that the identity matrix can be decomposed as
\[
I = \hat B {\hat B}^T + \sum_{i=1}^{n-1} {\hat B}_i {{\hat B}_i}^T
\] and therefore
\[
| I - {\hat B} {\hat B}^T | = | \sum_{i=1}^{n-1} {\hat B}_i {{\hat B}_i}^T | = 0
\]
which concludes that the eigenvalue $\lambda_i$ has been set properly.

\subsubsection{$Q_s^T$ are orthogonal generators of the affine hyperplane associated with $\lambda_i$}

Recalling the construction of $Q_s^T$ as the part of the qr  decomposition of
\[(A - \lambda_i I)^T Q \]
that is associated to nonzero values of the $R$ matrix,  
where $Q^T$ is an orthogonal annihilator of the $B$ vector, i.e. $Q^T B = 0$.
To show this, introduce the change of coordinates 
\[z = (A - \lambda_i I) x\]
It is well defined (its inverse is well defined) because $\lambda_i$ is not an eigenvalue of $A$. We claim that the hyperplane in these coordinates is generated by the complementary vectors $\hat B_i$, $i=1, \ldots n$ to $\hat B$, i.e. $\hat B ^T \hat B_i = 0$. This means that any gain vector assigning $\lambda_i$ can be expressed as (with arbitrary $\alpha_i \in \mathbb{R}, i=1,\ldots, n-1$)
\[
K = \left( \frac{\hat B}{\| B \|} + \sum_{i=1}^{n-1} \alpha_i \hat B_i \right)^T ( A - \lambda_i \, I)
\]
Indeed, by direct computation similar to the previous paragraph
\[
| \lambda_i A - A + B K | = | \lambda_i I - A + \hat B^T \hat B ( A - \lambda_i B) | = 0
\]
because $\hat B_i^T \hat B = 0$. This means that the hyperplane is generated in the orginal coordinates $x$ by 
\[
\text{span } \{ (A - \lambda_i I)^T \hat B_i | i=1,\ldots, n-1 \}
\]
which is $\text{span } Q_s^T$.

\subsubsection{The eigenvalues of the quotient are carried over to the the system prior to taking the quotient when using formula (\ref{pull})}

We must show that Formula (\ref{pull}) leaves unchanged the current well placed eigenvalue $\lambda_i$. The proof is by induction, so that assuming that $\bar K$ has been set properly in the sense that at the $i$-th step
\[
\lambda( \bar A - \bar B \bar K) = \{\lambda_n, \ldots \lambda_{n-i}\}
\]
then we must show that with $K_o$ the gain fixing $\lambda_{n-i-1}$ during the constructive step (and leaving the other eigenvalues yet to be defined), the following gain (changing the notation by droping the suffix 's', that is, $Q_s^T$ is from now on written $Q^T$ to simplify notation)
\[
K = K_o + \bar K Q^T
\]
assigns all the eigenvalues properly, with $A$ and $B$ designating the quotient at level one step above $\bar A$ and $\bar B$, i.e. of dimension $n-i-1$, i.e.
\[\lambda( A - B K) = \{\lambda_n, \ldots, \lambda_{n-i-1}\}
\]
By direct development we must show that
\[
\left| \lambda I - (A - B (\bar K Q^T + K_o)) \right| = \prod_{k=n-i-1}^n (\lambda - \lambda_k)
\]

\begin{eqnarray}
\left| \lambda I - (A - B (\bar K Q^T + K_o)) \right|  &=& \nonumber \\
\left| \left( \begin{array}{c} q^T \\ Q^T \end{array} \right) 
\left( \lambda I - A + B \bar K Q^T + B K_o \right) \left( \begin{array}{cc} q & Q \end{array} \right) \right| 
&=&  \nonumber \\
\left| \begin{array}{cc} \lambda - q^T A q + q^T B K_o q & - q^T A q + q^T B K_o Q \\
- Q^T A q + Q^T B K_o q & \lambda I - \bar A + \bar B \bar K + Q^T B K_o Q \end{array} \right| \label{bigDet} 
\end{eqnarray}
Now recalling the expression of $K_o$ and adjusting for the notations
\[
K_o = \frac{B^T}{B^T B} ( A - \lambda_{n-i-1} I) q q^T
\]
the top left scalar in the determinant (\ref{bigDet})  becomes
using the decomposition of the identity matrix
$I = \frac{B B^T}{B^T B} + \sum_{k=1}^{n-i-1} \hat B_k \hat B_k^T$ with
the property that $\hat B_k^T q = 0$
\begin{eqnarray}
\lambda - q^T A q + q^T B K_o q &=& \nonumber \\
 \lambda - q^T A q + q^T \frac{B B^T}{B^T B} (A - \lambda_{n-i-1} I) q &=&  \nonumber \\
\lambda - q^T A q + q^T (I - \sum_k \hat B_k \hat B_k^T) (A - \lambda_{n-i-1} I) q &=& \nonumber \\
\lambda - q^T A q + q ^T A q - \lambda_{n-i-1} q^T q = \lambda - \lambda_{n-i-1} \nonumber 
\end{eqnarray}
Another important result is that the first left column just below the first entry is zero in (\ref{bigDet}).
Indeed, since $\hat B_k^T q = 0$
\begin{eqnarray}
- Q^T A q + Q^T B K_o q &=& - Q^T A q + Q^T \frac{B B^T}{B^T B} (A - \lambda_{n-i-1} I) q q^T q \nonumber \\
&=& - Q^T A q + Q^T (I - \sum_k \hat B_k \hat B_k^T) (A - \lambda_{n-i-1} I) q \nonumber \\
&=& - Q^T A q + Q^T A q = 0 \nonumber
\end{eqnarray}
Finally, we show that 
\[
Q^T B K_o Q = Q^T \frac{B B^T}{B^T B} (A - \lambda_{n-i-1} I) q q^T Q = 0
\]
because $q^T Q = 0$.
This means that (\ref{bigDet}) is the product of the top left scalar times the determinant of the bottom part, that is (\ref{bigDet}) is equal to 
\[
(\lambda- \lambda_{n-i-1})\,\, | \lambda I - \bar A + \bar B \bar K | 
\]
completing the proof that Formula (\ref{pull}) leaves unchanged the current well placed eigenvalue during the descending phase of the quotient algorithm (with the current notation it is $\lambda_{n-i-1}$) and copying the eigenvalues that are placed using the quotient (those set by $\bar K$ on the quotient system $\bar A$, $\bar B$).

\subsubsection{Multiple identical real eigenvalues}

\subsubsection{Consequence of the loss of controllability on the algorithm}

\section{Algorithms based on obtaining the chain of matrices $C_{n-1} C_{n-2} \ldots C_1$}

\subsection{Chain of anchors in Ackermann's formula, the second algebroid method}

\label{algebroid2Algorithm}

The second algebroid method proceeds by taking a quotient to cancel the $B$ vector in the quotient space.
The quotients are indexed by $i$ ranging from $i = 0$ to $i = n-1$. The initial input vector is $B_0 = B$.
The input vector in the first quotient space is the image of $A_0 B_0$ through the first anchor, that is
$B_1 = \text{an}_1 A_0 B_0$. The anchor $\text{an}_1$ is chosen to cancel $B_0$ and to improve the conditioning
of the map $A_{t,1} = \text{an}_1 A_0$ with respect to $A_0$.The quotient $A$ matrix, which is $A_1$ (it will not be explicitly computed) induces  a mirroring effect of $A$ on the vectors modulo the $B$ vector and the vectors in the quotient operated by the $A_1$ matrix. Proceeding inductively the dimensions reduces and the size of the system as well. Taking the quotient introduces:
\begin{itemize}
\item  projection operators or anchors $\text{an}_i$
\item  quotient operators $A_{t,i}$
\item  input vectors in the quotient $B_{i}$
\end{itemize}
We will first give the algorithm and then provide the diagrams to explain the anchors $\text{an}_i$ and the significance of $A_{t,i}$.

\subsection{Description of the second algebroid algorithm}

\subsubsection{Forward sweep}
\vskip 0.5cm
\begin{itemize}
\item[{\bf Init}:] \begin{itemeqnarray}
A_0 &=& A \nonumber \\
B_0 &=& B \nonumber 
\end{itemeqnarray}
\item[{\bf Ind}:]
\begin{enumerate}
\item QR decomposition: \begin{itemeqnarray}
P\, R &:=& B_i \nonumber \\
\left(
\begin{array}{c}
q^T \\
Q^T
\end{array} \right) &:=& P \nonumber
\end{itemeqnarray}
$Q^T$ is of rank $n-1-i$ and cancels $B_i$, i.e. $Q^T \, B_i = 0$.
\item Scaling: \begin{itemeqnarray}
U \, \Sigma \, V^T &:=& Q_T \, A_t \nonumber \\
W &:=& U^T \nonumber
\end{itemeqnarray}
\item Storage:
\begin{itemeqnarray}
\text{an}_i &:=& W \, Q^T \nonumber \\
A_{t,i} &:=& \text{an}_i \, A_{t,i-1} \nonumber 
\end{itemeqnarray}
\item Quotients:
\begin{itemeqnarray}
A_{t,i} &:=& \text{an}_i \, A_{t,i-1} \, A \nonumber \\
B_{i} &:=& A_{t,i} \, B \nonumber 
\end{itemeqnarray}
\item Iteration: \begin{itemequation} i := i+1 \end{itemequation}
\end{enumerate}
\item[{\bf Ter}:] Stop when \begin{itemequation} i=n \end{itemequation} that is, when $A_{t,n-1}$ and $B_{n-1}$ are computed.
\end{itemize}
The forward sweep can be found as the code \texttt{buildOP} in Section \ref{buildOpCode}.

\subsubsection{The construction phase}
The constructive phase of the gain can either be done while computing the anchors and quotient operators,
or after that construction. Denoting the characteristic polynomial to be placed by
\[
P(\lambda) = \lambda^{n} + p_{n-1} \lambda^{n-1} + \ldots p_n
\] 
It proceeds along the following iterative steps

\begin{enumerate}
\item[{\bf Init}:]
\[K_0 := p_n \, I \qquad i=1 \]
\item[{\bf Iter}:]
\[K_{i} := A_{t,i}  p_{n-i-1} + \text{an}_i \, K_{i-1} \qquad i=i+1\]
\item[{\bf Ter}:] 
\begin{itemeqnarray}
\text{Stop when }i&=&n \nonumber \\
K_n &:=& K_{n-1} + A_{t,n-1} A \nonumber \\
K &:=& - \frac{1}{B_n} \, K_n \nonumber  
\end{itemeqnarray}
\end{enumerate}
The code is given in Section \ref{fbkCompExtractCode} and repeated here for convenience. 
\begin{verbatim}
n = length(B);
pp = poly(P);
pp = pp(end:-1:1); % put in reverse order
Kt = pp(1)*eye(n); % cumulative contribution to the final gain 
for i=1:n-1
Ot = Oo(((i-1)*n+1):((i-1)*n+(n-i)),1:n-i+1);
At = Ao(((i-1)*n+1):((i-1)*n+(n-i)),:);
Kt = At*pp(i+1) + Ot*Kt;
end;
Kt = Kt+At*A;
K = -1/(At*B)*Kt;
\end{verbatim}

\subsection{Diagrams for explaining the algorithm}

The quotients can be put into a diagram, where the maps $A_i$, $i=1,n-1$ represent the effect of the $A$ map on 
the quotient at the level $i$. We set $A_0 = A$.
\begin{center}
\begin{tikzcd}
\mathbb{R}^n \arrow{r}{A_0}  \arrow{d}{\text{an}_1}& \mathbb{R}^n \arrow{d}{\text{an}_1} \\
\mathbb{R}^{n-1} \arrow{r}{A_1} \arrow{d}{\text{an}_2} & \mathbb{R}^{n-1} \arrow{d}{\text{an}_2} \\
\mathbb{R}^{n-2} \arrow{r}{A_2} \arrow{d}{\text{an}_3} & \mathbb{R}^{n-2} \arrow{d}{\text{an}_3}\\
\vdots \arrow{d}{\text{an}_{n-1}} & \vdots \arrow{d}{\text{an}_{n-1}} \\
\mathbb{R} \arrow{r}{A_{n-1}} & \mathbb{R} 
\end{tikzcd}
\end{center}

When coding the algorithm, it is less computationally demanding to resort to the map from the original space
to the quotient space that renders commutative the effect of the original $A$ map on the quotient.

\begin{center}
\begin{tikzcd}
\mathbb{R}^n \arrow{r} {A_0} \arrow{rd}{A_{t,1}}  & \mathbb{R^n} \arrow{d}{\text{an}_1} \\
 & \mathbb{R}^{n-1} 
\end{tikzcd}
\end{center}

and $A_{t,k}$ corresponds to the map making the following diagram commutative

\begin{center}
\begin{tikzcd}
\mathbb{R}^n \arrow{r}{A_0^k} \arrow{rd}{A_{t,k}}  & \mathbb{R^n} \arrow{d}{\text{an}_k \circ \cdots \circ \text{an}_2 \circ \text{an}_1 } \\
 & \mathbb{R}^{n-k} 
\end{tikzcd}
\end{center}

Another advantage (apart from fewer operations) is that computing explicitly $A^k$ is avoided.
Associating the map 
\[
\text{an}_1 : \mathbb{R}^n \rightarrow \mathbb{R}^n 
\]
with the matrix 
\[
\text{an}_1
\]
and the composition of maps $\text{an}_2 \circ \text{an}_1$ with the matrix product $\text{an}_2 \, \text{an}_1$,
we have the following equalities
\begin{eqnarray}
\text{an}_1 \,  A_0   &=& A_{t,1}  \nonumber \\
\text{an}_2 \, \text{an}_1 \, A_0^2  &=& A_{t,2}  \nonumber \\
\text{an}_3 \, \text{an}_2 \, \text{an}_1 \, A_0^3  &=& A_{t,3}  \nonumber \\
&\vdots& \nonumber \\
\text{an}_k \, \cdots \, \text{an}_2 \, \text{an}_1 \, A_0^k  &=& A_{t,k}  \qquad k=1, \ldots, n-1 \nonumber
\end{eqnarray}
The anchors $\text{an}_k$ are chosen such that
\[
\text{an}_k \, \cdots \,  \text{an}_2 \, \text{an}_1 \, A^{k-1} B = 0 \qquad k=1, \ldots, n-1
\]
which gives the diagram of Figure \ref{anchorDiagram} where arrows are maps given by matrices and the nodes are specific vectors in the quotients.

\begin{figure}[ht]
\begin{center}
\begin{tikzcd}
B_0 \arrow{r}{A_0} \arrow{d}{\text{an}_1}  & A_0 B_0 \arrow{r}{A_0} \arrow{d}{\text{an}_1}& A_0^2 B_0 \arrow{r}{A_0} \arrow{d}{\text{an}_1} & \cdots \arrow{d}{\text{an}_1} \arrow{r}{A_0} & A_0^{n-1} B_0 \arrow{d}{\text{an}_1}\\
0 \arrow{r}{A_1}& B_1 \arrow{d}{\text{an}_2} \arrow{r}{A_1}& A_1 B_1 \arrow{r}{A_1} \arrow{d}{\text{an}_2}  & \cdots\arrow{d}{\text{an}_2}  \arrow{r}{A_1} & A_1^{n-2} B_1\arrow{d}{\text{an}_2}  \\
  & 0 \arrow{r}{A_2} & B_2 \arrow{d}{\text{an}_3} \arrow{r}{A_2}  & \cdots \arrow{d}{\text{an}_3}  \arrow{r}{A_2}  & A_2^{n-3} B_2 \arrow{d}{\text{an}_3}  \\
  &   & 0 \arrow{r}{A_3}  & \cdots \arrow{r}{A_3} \arrow{d}{}  & A_3^{n-4} B_3 \arrow{d}{} \\
  &   &   & \vdots \arrow{d}{\text{an}_{n-1}}   & \vdots \arrow{d}{\text{an}_{n-1}} \\
  &   &   & 0 \arrow{r}{A_{n-1}}  & B_{n-1} 
\end{tikzcd}
\end{center}
\caption{The diagram that links the chain of anchors $\text{an}_{n-1} \circ \ldots \circ \text{an}_2 \circ \text{an}_1$ to the factorization $C_{n-1} C_{n-2} \ldots C_1 = c$ of the last row of the inverse of the controllability matrix. \label{anchorDiagram}}
\end{figure}

\subsubsection{Controllability lemmas}
From the diagram of Figure \ref{anchorDiagram} it is straightforward to infer the following two lemmas. 

\begin{lemma}
\label{lemmaOne}
\begin{eqnarray*}
&\text{dim span} \{B_0, A_0 B_0, \cdots, A_0^{n-1} B_0 \} = n \nonumber  \\
&\Downarrow \nonumber \\
&\text{dim span} \{B_k, A_k B_k, \cdots, A_k^{n-k-1} B_k\} = n-k \qquad \forall k=1,\ldots,n-1 \nonumber
\end{eqnarray*}
\end{lemma}

The conclusion follows from the commutativity property of the diagram in Figure \ref{anchorDiagram}.
To make the proof precise, one can proceed by induction on $k$, and using the fact that $\text{an}_k$ is
of same rank as $Q_k^T$ which is full row rank $n-k-1$ and cancels $B_k$. Indeed, if one supposes
the induction hypothesis
\[
\text{dim span} \{ B_k, A_k B_k, \cdots, A_k^{n-k-1} B_k \} = n-k
\]
then 
\begin{equation}
\label{collection}
\text{dim span} \{ Q_k A_k B_k, \cdots Q_k A_k^{n-k-1} B_k\} = n-k-1
\end{equation}
since $Q_k$ is of rank $n-k-1$ and cancels vectors proportional to $B_k$, and $B_k$ has been left out in (\ref{collection})---
furthermore 
 $A_k B_k$, $A_k^2 B_k$, $\ldots$, $A_k^{n-k-1} B_k$ are linearly independent by the induction hypothesis. Then with $\text{an}_k = W_k \, Q_k^T $, the fact that $B_{k+1} = \text{an}_k \, A_k B_k$, and 
the commutativity of the diagram between row $k$ and row $k+1$
\[
\text{dim span} \{ B_{k+1}, A_{k+1} B_{k+1}, \cdots A_{k+1}^{n-2-k} B_{k+1} \} = n-k-1
\]
as well, proving the induction step and the lemma (the initialisation of the induction argument is the same, just set $k=0$).

The contrapositive statement  
to Lemma \ref{lemmaOne} is useful for testing the loss of  controllability:

\begin{lemma}
\begin{eqnarray}
\exists k \neq 0, B_k = 0 \Rightarrow \text{rank} \left( \begin{array}{cccc} B_0 & A_0 B_0 & \cdots 
A_0^{n-1} B_0 \end{array} \right) < n \nonumber
\end{eqnarray}
\end{lemma}
The proof is straightforward since $B_k = 0$ prevents 
achieving 
\[\text{dim span} \{B_k, A_k B_k, \cdots, A_k^{n-k-1} B_k\} = n-k\]
 and the contrapositive of Lemma \ref{lemmaOne} leads to the loss of controllability.
A similar controllability test appeared in \cite{MullhauptCDC}. 

\subsection{The theorem of the second algebroid algorithm and its corollary}

A direct consequence of the diagram of Figure \ref{anchorDiagram} gives a chain of matrices giving the last row 
of the controllability matrix.

\begin{theorem}
\label{algebroid2Theorem}
Le $c = e_n^T \mathcal{C}^{-1}$ be the last row of the inverse of the controllability matrix. A factorization
according to $c = C_{n-1} C_{n-2} ... C_1$ mentioned in Section \ref{AckermannFactorization} is obtained from the diagram of Figure \ref{anchorDiagram}, and is a result of applying the second algebroid algorithm:
\begin{equation}
\label{cFormula}
c = \frac{1}{B_{n-1}} \, \text{an}_{n-1} \, \text{an}_{n-2} \, \ldots \, \text{an}_2 \, \text{an}_1 
\end{equation}
\end{theorem}
To prove this statement, notice that $c \, A_0^{n-1} B_0 = 1$ because of the definition of $c$. The last column of the diagram in Figure \ref{anchorDiagram} gives
\[ B_{n-1} = \text{an}_{n-1} \ldots \text{an}_2 \text{an}_1 \, A_0^{n-1} \, B_0 \]
which is after dividing both sides by $B_{n-1}$ the relation $c \, A_0^{n-1} B_0 = 1$ if $c$ equals the value given by  (\ref{cFormula}).

\begin{corollary}
The feedback gain $K$ can be computed using the nested sequence of Section \ref{chainAckermannFactorization}
\end{corollary}

This is a direct consequence of the chain of anchors in (\ref{cFormula}) which is the same as the same--modulo the factor $1/B_{n-1}$--as the chain of matrices $C_{n_1} C_{n_2} \ldots C_1$ of Section \ref{chainAckermannFactorization}. 

Taking into account the factor $1/B_{n-1}$ the procedure of Section \ref{chainAckermannFactorization} leads to the constructive phase of the second algebroid algorithm. With this last theorem, the correctness of the second algebroid algorithm has been proved.

\subsection{Example to illustrate the second algebroid algorithm}

Applying the forward sweep part of the algorithm coded in  \texttt{buildOp} (see Section \ref{buildOpCode}) on our "fil rouge" example produces:
\begin{verbatim}
> A0
A0 =
     1     3     5
     7    13    17
     1     1     1
> B0
B0 =
     1
     1
     1
> Oo
Oo =
    0.2673   -0.8018    0.5345
   -0.7715    0.1543    0.6172
    0.0000    0.0000    0.0000
   -0.0240   -0.9997    0.0000
    0.0000    0.0000    0.0000
    0.0000    0.0000    0.0000
> Ao
Ao =
   -4.8107   -9.0869  -11.7595
    0.9258    0.3086   -0.6172
    0.0000    0.0000    0.0000
   -0.5399   -2.6995   -4.6792
    0.0000    0.0000    0.0000
    0.0000    0.0000    0.0000
> 
\end{verbatim}
The anchors $\text{an}_i$, $i=1,2$ are stored in the \texttt{Oo} array. The maps $A_{t,i}$, $i=1,2$ are stored in the \texttt{Ao} array. Let us illustrate the results obtained by cross-checking a few identities.
The identity
\begin{eqnarray}
\text{an}_1 \, A_0 = A_{t,1} \nonumber 
\end{eqnarray}
comes down to checking the equality between the following two statements
\begin{verbatim}
> Oo(1:2,:)*A0
ans =
   -4.8107   -9.0869  -11.7595
    0.9258    0.3086   -0.6172
> Ao(1:2,:)
ans =
   -4.8107   -9.0869  -11.7595
    0.9258    0.3086   -0.6172
> 
\end{verbatim}  
These expressions match. The following identity 
\[
\text{an}_2 \, \text{an}_1 \, A_0^2 = A_{t,2}  
\]
is also checked as 
\begin{verbatim}
> Oo(4,1:2)*Oo(1:2,:)*A0*A0
ans =
   -0.5399   -2.6995   -4.6792
> Ao(4,:)
ans =
   -0.5399   -2.6995   -4.6792
\end{verbatim}
match as well, and finally examining the diagram \ref{anchorDiagram}, we cross-check
\[
B_1 = \text{an}_1 \, A_0 \,  B_0 = A_{t,1} \, B_0
\]
\begin{verbatim}
> Ao(1:2,:)*B0
ans =
  -25.6571
    0.6172
> Oo(1:2,:)*A0*B0
ans =
  -25.6571
    0.6172
\end{verbatim}
and
\[
B_2 = \text{an}_2 \, \text{an}_1 A_0^2 \, B_0 = A_{t,2} \, B_0
\]
\begin{verbatim}
> Ao(4,:)*B0
ans =
   -7.9186
> Oo(4,1:2)*Oo(1:2,:)*A0*A0*B0
ans =
   -7.9186
\end{verbatim}
The second part of the algorithm--the constructive gain phase--produces 
\begin{verbatim}
% the desired characteristic polynomial coefficients in reverse order
pp =
     6    11     6     1
% The initialisation stage:
Kt =
     6     0     0
     0     6     0
     0     0     6
% the iterations Kt = At*pp(i+1) + Ot*Kt give:
Kt =
  -51.3142 -104.7664 -126.1473
    5.5549    4.3205   -3.0861
Kt =
   -7.5587  -17.9969  -21.9562
% and the termination stage:
> Kt = Kt+At*A
Kt = 
  -31.6745  -59.3896  -75.2268
> K = -1/(At*B)*Kt
K = 
   -4.0000   -7.5000   -9.5000
\end{verbatim}
and we obtain the correct eigenvalues for $A + B K$.
\begin{verbatim}
> eig(A+B*K)
ans =
   -3.0000
   -2.0000
   -1.0000
\end{verbatim}

\section{Numerical experiments}

\subsection{Integer number example}

The first example is an example that rapidly leads to some numerical problems in representing the final gain
vector properly. This comes from the position of the eigenvalues and the rapidly ill conditioned controllability matrix.
Geometrically, using the understanding from the affine hyperplanes, as the dimension increases, the intersection 
of the affine hyperplanes is situated on a point that requires both high values for some of their components and small values for other components. 

The examples are generated by the following lines of code
\begin{verbatim}
 NN = 10;
 AA = [1:NN; [eye(NN-1),ones(NN-1,1)]];
 AA(3:end,1) = -ones(NN-2,1);
 BB = ones(NN,1);
\end{verbatim}
which corresponds to 
\[
A = \left(\begin{array}{cccccc} 
1 & 2 & 3 & \ldots & n-1 & n \\
1 & 0 & 0 &\ldots & 0 & 1 \\
-1 & 1 & 0 &\ldots & 0 & 1 \\
-1 & 0 & 1 & 0 & \ldots & 1 \\
\vdots & 0 & \ddots & \ddots & \ldots & 1 \\
-1 & 0 & 0 & \ldots & 1 & 1 
\end{array}
\right) \qquad B = \left( \begin{array}{c} 1 \\ 1 \\ 1 \\ 1 \\ \vdots \\ 1   \end{array} \right) 
\]
The desired eigenvalues are $\lambda_1 = -1$, $\lambda_2 = -2$, $\ldots$, $\lambda_n = -n$ in that order.
Using 64 bits arithmetic with $\epsilon = 2.2204 \, 10^{-16}$ the following results of algebroid algorithm 1 ('alg1', projection
into the hyperplanes) and algebroid algorithm 2 ('alg2', successive cancelation of the input vector $B_i$ of the quotient)  are produced in Sysquake and compared to Matlab's 'place' command. The resulting eigenvalues are represented in the following table.

\begin{center}
\begin{tabular}{c | c | c}
alg1 & alg2 & place \\
\hline
-9.999988096777126 & -9.999987985525921 & -10.000164107737568 \\
-9.000043978895723 & -9.000044227665359 & -8.999371863648731 \\
-7.999935955089022 & -7.999935953734287 & -8.000964166915645 \\
-7.000046221237703 & -7.000045767406014 & -6.999242490281093  \\
-5.999983141210241 & -5.999983658058601 & -6.000322564954156 \\
-5.000002682351759 & -5.000002430833500 & -4.999928035624746 \\
-3.999999942133777 & -3.999999999694501 & -4.000007539812596 \\
-2.999999981610343 & -2.999999976182774 & -2.999999713947878 \\
-2.000000000712693 & -2.000000000847365 & -2.000000002108917 \\
-0.999999999998041 & -0.999999999997960 & -0.999999999951879
\end{tabular}
\end{center}
Although there is some inaccuracy, the algorithms are comparable, and the pole placement can be
judged satisfactory.
Increasing the dimension by one unit setting $n = 11$ with the instruction
\begin{verbatim}
NN = 11;
\end{verbatim}
leads to the the following tables
\begin{center}
\begin{tabular}{c | c}
alg1 & alg2 \\
\hline
   -11.009927094174800 &  -10.846862760359385 \\
    -9.957181536060149 &  -10.384497587913568  \\
    -9.067188389254070 &  -8.468508021828071 + 0.498862941311272j \\
    -7.944533081061675 &  -8.468508021828071 - 0.498862941311272j \\
    -7.026358424804697 &  -6.760343168864634 \\
    -5.994389988439849 & -6.076040368772530 \\
    -5.000380816161343 &  -4.995792337226679  \\
    -4.000046379477965 &  -3.999367994633918 \\
    -2.999994187410914 & -3.000081315485191  \\
    -2.000000103128592 &  -1.999998421821761 \\
    -0.999999999928366 &  -1.000000001278820    
\end{tabular}
\end{center}

\begin{center}
\begin{tabular}{c}
place \\
\hline
 -11.121510720787899 + 0.000000000000000i \\
  -9.620049837351161 + 0.579008846167512i \\
  -9.620049837351161 - 0.579008846167512i \\
  -7.375497284454621 + 0.511289713821731i \\
  -7.375497284454621 - 0.511289713821731i \\
  -5.860721968875364 + 0.000000000000000i \\
  -5.028671908237159 + 0.000000000000000i \\
  -3.997965809156488 + 0.000000000000000i \\
  -3.000056435717809 + 0.000000000000000i \\
  -1.999999711062670 + 0.000000000000000i \\
  -0.999999998825601 + 0.000000000000000i \\
\end{tabular}
\end{center}

The error is significantly more pronounced, and the pole placement can be judged as borderline. Some of the eigenvalues have imaginary parts. There is no such bifurcation for 'alg1' there is only one bifurcation for 'alg2'. There are two such bifurcations when applying 'place' of Matlab. Increasing to $n=12$ gives three bifurcations for 'alg1' and 'alg2' and 'place' has 4 such bifurcations. The pole placement is unsatifactory for all algorithms.

\begin{center}
\begin{tabular}{c|c|c}
alg1 & alg2 & place \\
\hline
 -12.7034 & -12.7137 &   -12.7647 + 1.5577i \\
 -10.9301 + 1.7027j & -10.9338 + 1.7163j &   -12.7647 - 1.5577i \\
 -10.9301 - 1.7027j & -10.9338 - 1.7163j &  -9.4596 + 2.7072i \\ 
 -8.0919 + 1.5917j & -8.0874 + 1.6025j &   -9.4596 - 2.7072 \\
 -8.0919 - 1.5917j & -8.0874 + 1.6025j &   -6.7149 + 1.7711i \\
 -6.1326 + 0.4986j & -6.1281 + 0.5032j &  -6.7149 - 1.7711i \\
 -6.1326 - 0.4986j & -6.1281 - 0.5032j &   -5.0730 + 0.5385i \\
 -4.9911 & -4.9916 &    -5.0730 - 0.5385i \\
 -3.9961 & -3.9959 &   -3.9703 + 0.0000i \\
 -3.0003 & -3.0003 &   -3.0006 + 0.0000i \\
 -2.0000 & -2.0000 &   -2.0000 + 0.0000i \\
 -1.0000 & -1.0000 &   -1.0000 + 0.0000i \\
\end{tabular}
\end{center}

Placing the eigenvalues in reverse order $\lambda_1 = -n$, $\lambda_2 = -(n-1)$, $\ldots$, $\lambda_{n-1} = -2$, $\lambda_{2} = -1$ with the command
\begin{verbatim}
-(NN:-1:1)
\end{verbatim}
changes the results produced by algorithms 'alg1' and 'place' of Matlab. The algorithm 'alg2' using the characteristic polynomial without its factorization is insensitive to the change of order. For example, for the last case with $n=12$, 'place' is improved because there are 3 bifurcations instead of 4, and 'alg1' remains with 3 bifurcations.
\begin{center}
\begin{tabular}{c|c|c}
alg1 & alg2 &  place \\
\hline
  -12.7033 &-12.7137  & -13.8776 + 0.0000i \\
-10.9296 + 1.7025j &  -10.9338 + 1.7163j &   -11.3620 + 2.9956i \\
 -10.9296 - 1.7025j & -10.9338 - 1.7163j & -11.3620 - 2.9956i \\
  -8.0912 + 1.5904j & -8.0874 + 1.6025j & -7.7667 + 2.6993i \\
   -8.0912 - 1.5904j & -8.0874 - 1.6025j & -7.7667 - 2.6993i \\
    -6.1332 + 0.4945j & -6.1281 + 0.5032j & -5.6377 + 1.3136i \\
    -6.1332 - 0.4945j & -6.1281 - 0.5032j &-5.6377 - 1.3136i \\
    -4.9925 & -4.9916  & -4.4996 + 0.0000i \\
    -3.9959  & -3.9959  & -4.0870 + 0.0000i \\
    -3.0003 & -3.0003 & -2.9989 + 0.0000i \\
   -2.0000  & -2.0000  & -2.0000 + 0.0000i \\
    -1.0000 & -1.0000 & -1.0000 + 0.0000i 
\end{tabular}
\end{center}

Similar sensitivities in the order of the eigenvalues  occur for 'alg1' and 'place' for  $n=10$ and $n=11$. Interestingly, 'alg1' that seemed to outperform the two others when $n=11$---because it did not display a single bifurcation---does bifurcate when the eigenvalues are presented in 
reverse order $-11$, $-10$, $\ldots$, $-1$:
\begin{center}
\begin{tabular}{c}
alg1 \\
\hline
 -11.0611         \\
 -9.5638 + 0.2258j \\
 -9.5638 - 0.2258j \\
 -7.5598        \\
 -7.2878        \\ 
 -5.9604        \\ 
 -5.0029        \\ 
 -4.0004        \\ 
 -2.9999        \\ 
 -2.0000        \\
 -1.0000
\end{tabular}
\end{center}

\subsection{Ring arithmetic results}

Since all entries in $A$ and $B$ are integers and the eigenvalues are integers as well, the previous example can serve as a test for the ring arithmetic version of the 'alg2'. 

An interesting property of the ring arithmetic version of 'alg2' in such a ring arithmetic setting is that it computes an exact solution to the problem. It can then be used as a 'ground truth' solution and identify whether the inaccuracies are intrinsic due to the algorithm used to place the eigenvalues or are due to either inherent limitation of the gain and/or the algorithm to compute the eigenvalues. If the errors are solely from the algorithm to place the eigenvalues, then the ring arithmetic version should improve the observed limit of $n=11$. If it is inherent to the impossibility to represent sucessfully the entries of the gain vectors (due to the arithmetic limitation of 64 bits precision) and/or the algorithm to compute the eigenvalues numerically, then the ring arithmetic version should not improve the limit observed of $n=11$. 

\subsubsection{$n = 11$, ring arithmetic}
Using the ring arithmetic version of 'alg2'
\begin{verbatim}
KK11 = [7817883664811469804057/297365203664055278341,...
 66347135266209260491107/297365203664055278341,...
 715307440643594285832987/297365203664055278341,...
 5108463570029711309325053/297365203664055278341,...
 24279372098464306568093845/297365203664055278341,...
 74798168434160582892384569/297365203664055278341,...
 136845070738935394124936213/297365203664055278341,...
 106617412978197400238773250/297365203664055278341,...
 -(69104192347823610988017594/297365203664055278341),...
 -(186582984738415277335631860/297365203664055278341),...
 -(92730562359273966067064439/297365203664055278341)];
\end{verbatim}
produced the eigenvalues in the following table

\begin{center}
\begin{tabular}{c}
ring alg2 \\
\hline
  -11.046830477565503 \\
   -9.673388926331968 \\
   -9.429782589303520 \\
   -7.699694089022847 \\
   -7.178111035258130 \\
   -5.969756800127997 \\
   -5.002156766157062 \\
   -4.000321026021105 \\
   -2.999957426885713 \\
   -2.000000863838766 \\
   -0.999999999253550
\end{tabular}
\end{center}

The precision is not very accurate despite the exactness of the computation of \texttt{KK11}.
The result is still somewhat promising since there is no bifurcation. 

\subsubsection{$n = 12$, ring arithmetic}
However, increasing the dimension by one, the ring version of 'alg2' produces the exact value of the gain
\begin{verbatim}
KK12 = [3140867001984180016036461/100701343380251789934337,...
 32463700215024014546326491/100701343380251789934337,...
 433968633546560213091669147/100701343380251789934337,...
 3931398036873040592316764237/100701343380251789934337,...
 24528600373899823370244217765/100701343380251789934337,...
 104772649587412878088636414193/100701343380251789934337,...
 295598922877646668386365328773/100701343380251789934337,...
 499124346841391853303086344214/100701343380251789934337,...
 344789964075341274989916614646/100701343380251789934337,...
 -(290515578148790898307469121652/100701343380251789934337),...
 -(665350044862049195830462375466/100701343380251789934337),...
 -(317341775875018592857093471849/100701343380251789934337)];
\end{verbatim}
There is some disapointment, because there is no improvement on the floating point arithmetic versions of the pole placement algorithms. The following table gives the results of the eigenvalues computed with the command
\begin{verbatim}
eig(AA12-BB12*KK12)
\end{verbatim}

\begin{center}
\begin{tabular}{c}
ring alg2 \\
\hline
 -12.555248633420411  \\
 -10.867280335425483 + 1.494979442980424j \\
 -10.867280335425483 - 1.494979442980424j \\
 -8.148495215966406 + 1.404942415254940j \\
 -8.148495215966406 - 1.404942415254940j \\
 -6.209294613248384 + 0.365554425915730j \\
 -6.209294613248384 - 0.365554425915730j \\
 -4.997166798363346                    \\
 -3.997277647108256                   \\ 
 -3.000168276711908                    \\
 -1.999998316475768                    \\
 -1.000000000583342  
\end{tabular}
\end{center}

with the typical 3 bifurcations observed in the floating point arithmetic versions of the algorithms.

There is still one question remaining. Does the inacuracy observed--- despite the exact computation of the gain--- come from the inherent incapacity to represent the gain properly using 64 bits floating point arithmetic  during the process of computing the matrix $A - B K$--- which could potentially have severe rounding offs and cancelling--- or is it due to the computation algorithm for obtaining the  eigenvalues?

To try to settle this question, we will resort to a lower precision arithmetic taking a 32 bits version of Sysquake with $\epsilon = 1.1921 \, 10^{-7}$. The separation occurs between $n = 7$ and $n=8$ in this case. Running 'alg1' gives 

\begin{center}
\begin{tabular}{c|c}
alg1 & alg2 \\
\hline
 -8.013718 + 1.077389j &  -8.816705 \\
 -8.013718 - 1.077389j & -6.491629 + 1.586151j \\
 -5.060093 + 1.086392j & -6.491629 - 1.586151j \\
 -5.060093 - 1.086392j & -4.057369 + 0.406279j \\
 -3.922629 &  -4.057369 - 0.406279j \\
 -2.919564 &  -3.096357 \\
 -2.010736 &  -1.988744 \\
 -0.999894 &  -1.988744 
\end{tabular}
\end{center}

The first column has been achieved using the gain vector
{\tiny
\begin{verbatim}
eig(AA+BB*Ka)
Ka =
  -14.179460  -56.722354 -256.553253 -651.860656 -753.698608  105.431648 1003.696289  585.885925
\end{verbatim}}
(Attention to the sign, it is $A+B K$ not $A-B K$). To check the arithmetic roundoff, compute
{ \tiny
\begin{verbatim}
> AA+BB*Ka
ans =
  -13.179460  -54.722354 -253.553253 -647.860656 -748.698608  111.431648 1010.696289  593.885925
  -13.179460  -56.722354 -256.553253 -651.860656 -753.698608  105.431648 1003.696289  586.885925
  -15.179460  -55.722354 -256.553253 -651.860656 -753.698608  105.431648 1003.696289  586.885925
  -15.179460  -56.722354 -255.553253 -651.860656 -753.698608  105.431648 1003.696289  586.885925
  -15.179460  -56.722354 -256.553253 -650.860656 -753.698608  105.431648 1003.696289  586.885925
  -15.179460  -56.722354 -256.553253 -651.860656 -752.698608  105.431648 1003.696289  586.885925
  -15.179460  -56.722354 -256.553253 -651.860656 -753.698608  106.431648 1003.696289  586.885925
  -15.179460  -56.722354 -256.553253 -651.860656 -753.698608  105.431648 1004.696289  586.885925
\end{verbatim}}
and insert both results in Sysquake 64 bits precision, that is \texttt{Ka} and the result of
\texttt{AA+BB*Ka}, and compare the eigenvalues of the last expression with  the eigenvalues of \texttt{AA+BB*Ka} computed with 64 bits arithmetic, but keeping the \texttt{Ka} computed with 32 bits arithmetic. This leads to the following two columns in the table

\begin{center}
\begin{tabular}{c|c}
\texttt{AA+BB*Ka} 32 bits, \texttt{eig} 64 bits & \texttt{AA+BB*Ka} 64 bits, \texttt{eig} 64 bits \\
\hline 
 -8.162870655036357   &  -8.162870655036357  \\
  -6.530537765151010 + 0.529298100457248j & -6.530537765151010 + 0.529298100457248j \\
  -6.530537765151010 - 0.529298100457248j & -6.530537765151010 - 0.529298100457248j \\
  -4.633687310905854 & -4.633687310905854  \\
  -4.154158715594606  & -4.154158715594606 \\
  -2.988437062048782  & -2.988437062048782 \\
  -2.000239264801626  & -2.000239264801626  \\                  
  -1.000000461310615   & -1.000000461310615 
\end{tabular}
\end{center}
and there is only very small differences between the two columns.
Now let us use the  exact value of the gain 
\begin{verbatim}
KK8 = [519515210277/36638795621,...
 2078221618718/36638795621,...
 9399790968804/36638795621,...
 23883421055437/36638795621,...
 27614625334253/36638795621,...
 -(3862903459832/36638795621),...
 -(36774234975734/36638795621),...
 -(21466161518325/36638795621)];
\end{verbatim}
---(attention: this time it is $A - B K$, the minus sign)---
produces
{\tiny
\begin{verbatim}
KK8 =
   14.179373   56.721885  256.552947  651.861511  753.699035 -105.432037 -1003.696472 -585.886047
\end{verbatim}}
instead of \texttt{Ka} computed with 'alg1', but still stored in 32 bits arithmetic. After repeating the computation of the previous table, only changing \texttt{Ka} but computed from the 
exact value using 32 bits arithmetic, the same \texttt{Ka} for both first columns. The last column uses the exact solution \texttt{KK8} evaluated this time using 64 bits arithmetic and all expressions  is in 64 bits.
\begin{center}
\begin{tabular}{c|c|c}
\texttt{AA+BB*Ka} 32 bits, \texttt{eig} 64 bits & \texttt{AA+BB*Ka} 64 bits, \texttt{eig} 64 bits& \texttt{A+BB*KaExact} 64 bits, \texttt{eig} 64 bits  \\
 -9.424880 & -8.077604712315544 & -8.000000000255488 \\
 -6.486456 + 2.264287j & -6.680323904081507 & -6.999999999312520 \\
 -6.486456 - 2.264287j &   -6.322578204304940 & -6.000000000613400 \\
 -3.646944 + 0.711540j & -4.89637640832420 &  -4.999999999841751 \\
 -3.646944 - 0.711540j &  -4.025088703450608 &  -3.999999999954670 \\
 -3.349041 & -2.998196982904853 &  -3.000000000024483 \\
 -1.959010 &   -2.000025608029641 &    -1.999999999997695 \\
 -1.000558 &  -1.000000476588480 &  -1.000000000000016
\end{tabular}
\end{center}
The first column indicates that the eigenvalue computation is not accurate because setting $A - B K$ computed with the matrices  and vectors using 32 bits accuracy arithmetic  gives completely different values when computed accurately using 64 bits arithmetics (second column). However, there is still significant truncation errors that are reflected by the difference between  the last
two columns.

\subsection{Dynamical simulation results}

\label{testSection}
\begin{eqnarray}
\bar A = \left( \begin{array}{ccccc} 
1 & 0 & \ldots & 0 & 0 \\
0 & 2^{-2} & 0 &\ldots & 0 \\
\vdots & &\ddots & & \vdots \\
0 & \ldots & 0 & (n-1)^{-2} & 0 \\
0 & \ldots & 0 & 0 & n^{-2} 
\end{array}
 \right) \qquad \qquad \bar B = \left( \begin{array}{c}
1 \\ 1 \\ \vdots \\ 1 \\1
\end{array} \right)
\nonumber 
\end{eqnarray}

This is the structural matrices for the family of systems that will be generated through a random choice of similarity transform matrices that are chosen to be orthogonal. That is, let us asssume that 
\begin{eqnarray}
A = Q^T \bar A Q \qquad B = Q^T \bar B 
\end{eqnarray}
with $Q$ a randomly chosen orthogonal matrix (i.e. $Q^T Q = I$). 

The family of systems as $n$ ranges over the integers leads quickly to an ill-conditioned controllability matrix $\mathcal C$ when finite precision arithmetic is used (be it floating point or fixed point). Let us admit nevertheless that a row vector $c^T$ can be exactly computed that cancels the theoretical controllability matrix using infinite precision arithmetic. Clearly, the structure is controllable and it is only the numerical representation of the system that leads to the loss of controllability. 

Ackermann's formula cannot be applied beyond a given index $n=N$ which depends on the precision chosen for the arithmetic, even when the vector $c^T$ (that cancels the first $n-1$ columns of $\mathcal C$ and such that $c^T A^{n-1}B =1$) is exactly computed. An explanation for this is that Ackermann's formula necessitates the computations of high powers of the $A$ matrix. 

\section{The Varga algorithm comparison}

The complete algorithm is given in code in Section \ref{VargaAlgo}. 
The system is first brought into the Schur form. This form has the following interpretation when looked rom the perspective of the intersection of hyperplanes. Recall that all eigenvalues to be placed  are considered real and distinct and different from those of the initial matrix (before the Schur reduction which keeps the eigenvalues unchanged). 

\begin{itemize}
\item The Schur form $A_s$ of the system matrix $A$ has the particularity that an hyperplane appears coded in the last column of $A_s$. The associated last component of the $B_s$ vector can be used to translate the hyperplane to its value corresponding the eigenvalue to be assigned. 

\item Changing coordinates, by a suitable transformation,  one can bring another column of the transformed matrix without affecting the hyperplane that has been translated (and hence without affecting the eigenvalue corresponding to that hyperplane). 

\item The operations are repeated (translation of the hyperplane in the last column -- changing the eigenvalue -- then coordinate change to permute) until all eigenvalues are placed.

\item Bookkeeping the permutations and the scalar translations provides the feedback gain.
\end{itemize}

Hence, one sees that the idea of quotient is embedded in the successive transformation of the Schur forms while preserving the Schur structure. In a certain sense, the hyperplane in the last column can serve as the representative of the class of all those translated hyperplanes corresponding to a single eigenvalue. The shape of the hyperplane does not change but only the value of the eigenvalue by translating using the last component of the $B_s$ vector. Then through permutatation (change of coordinates) another class corresponding to another hyperplane is treated.

After this sketch of the interpretation of the Schur algorithm  from the point of view of the hyperplane intersection, we will provide the numerical experiments associated to the ones we have already performed with the comparison between the Miminis-Paige algorithm and the algorithms contributed in this paper (the two algebroid methods).

\section{Numerical simulation}

\subsection{The case study}

Using the $A$ and $B$ matrices described in Section \ref{testSection}, a linear time-invariant system is simulated using the classical Runge-Kutta algorithm
\[
\dot x = A x + B u
\]
with a given intial condition $x(0) = x_0$ and using either the $K$ vector $u = - K x$ or the feedback function $u = k(x)$.

Instead of producing a gain vector $K$ such that $A - B K$ has the suitable eigenvalues or desired characteristic polynomial, the algorithm produces a numerical algorithmic feedback function $k(x): \mathbb{R}^n \rightarrow \mathbb{R}$ accurate to machine precision (as an independent entity) without factorizing the function $k(x)$ as $- K x$. This means that the feedback $u = k(x)$ has to be implemented instead of $u = - K x$. It is more costly than a single row vector multiplication. Hence the algorithm has two stages, namely (i) a constructive off-line stage that builds suitable real orthogonal projection matrices and (ii) the feedback stage comprising the construction of the function $k(x)$ that should be implemented for the control scheme to be effective on-line.

The purpose is to check asymptotic stability and regularity of the solutions as the dimension $n$ of the matrices are increased. Depending on the algorithm used and the size $n$ instability and irregularities in the solution might appear.

In the following simulation results one notices large transients befor convergence. The transients augment as the size $n$ increases. Beyond a certain value of $n$ instability occurs and/or inacurate solutions appear.

\subsection{Without random orthogonal similarity transforms}

\begin{itemize}
\item All algorithms perform satisfactorily well in this section
\item Numerical instability occurs beyond $n=16$.
\end{itemize}

\subsection{With random orthogonal similarity transforms}

\begin{itemize}
\item Numerical instability is comparitively much worse in this scenario, and only the Miminis-Paige gives satisfacory results beyond $n=10$.
\item The Miminis-Paige gives an error starting to be significant for $n=15$ and that closed-loop instability occurs for $n=16$.
\item The second proposed algorithm gives the best error up to and including $n=17$. Strange noise appears for $n=18$ but maintaining closed-loop stability.
\end{itemize}

\begin{figure}
\includegraphics[width=6cm]{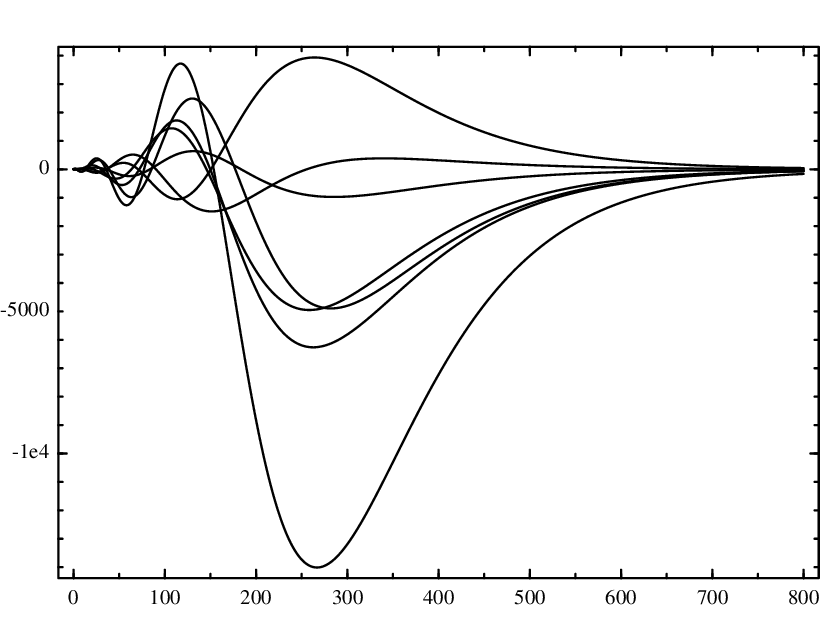}
\includegraphics[width=6cm]{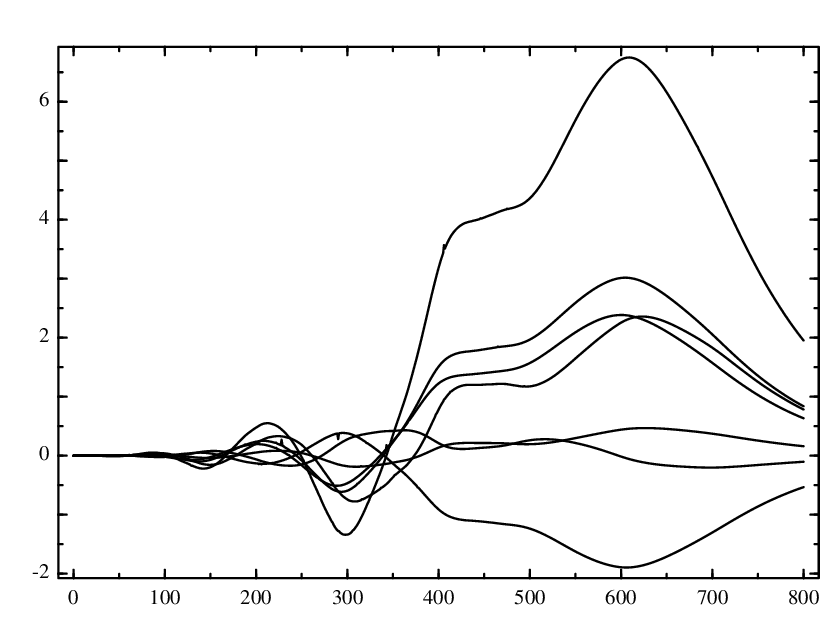}
\caption{$n=7$, left the state trajectory, all computations performed using 32 bits precision $\epsilon = 1.1921e-7$. Right, the error with 64 bits precision arithmetic $\epsilon = 2.2204e-16$}
\end{figure}

\begin{figure}
\includegraphics[width=6cm]{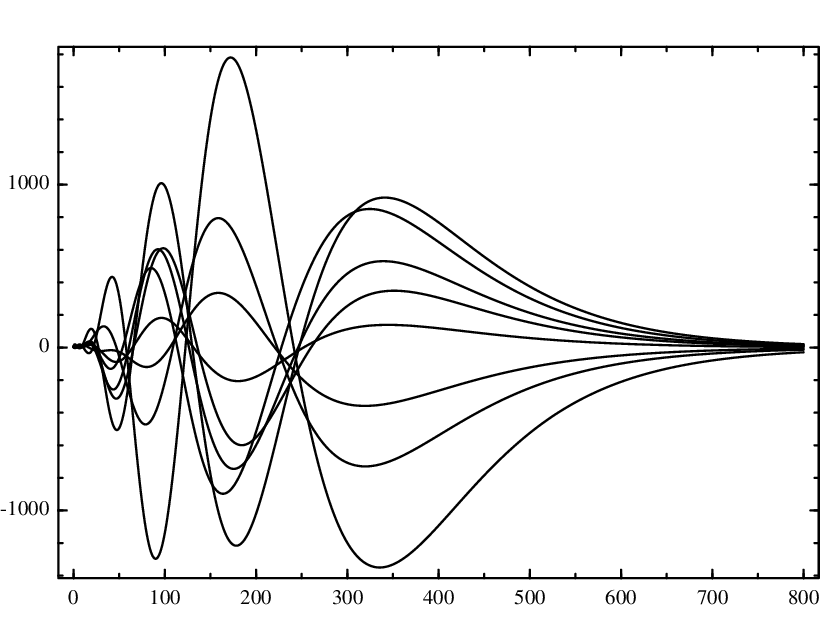}
\includegraphics[width=6cm]{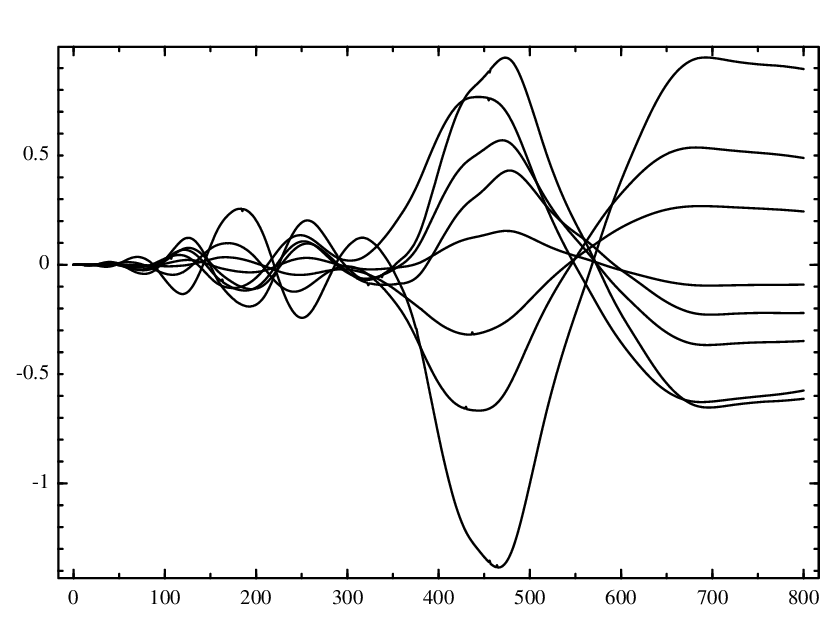}
\caption{$n=8$, left the state trajectory, all computations performed using 32 bits precision $\epsilon = 1.1921e-7$. Right, the error with 64 bits precision arithmetic $\epsilon = 2.2204e-16$}
\end{figure}

\begin{figure}
\includegraphics[width=6cm]{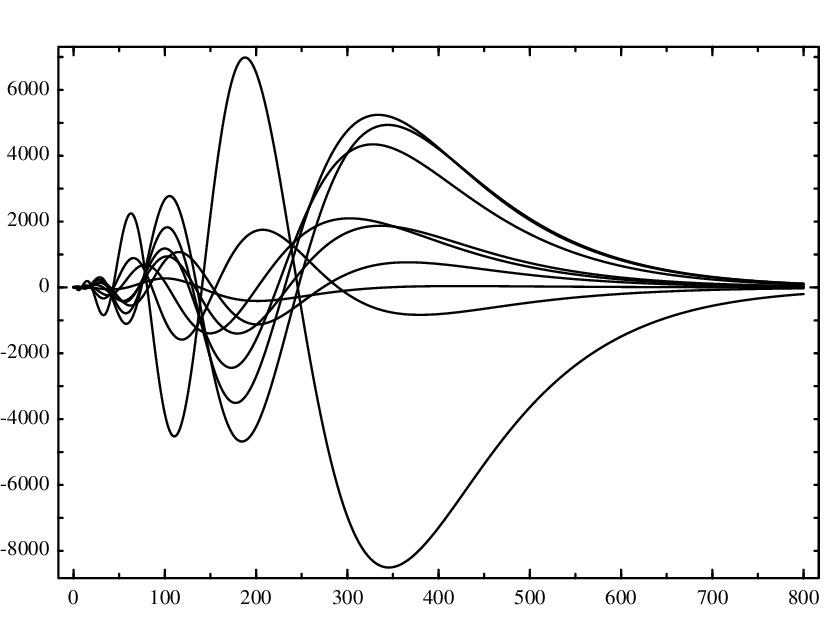}
\includegraphics[width=6cm]{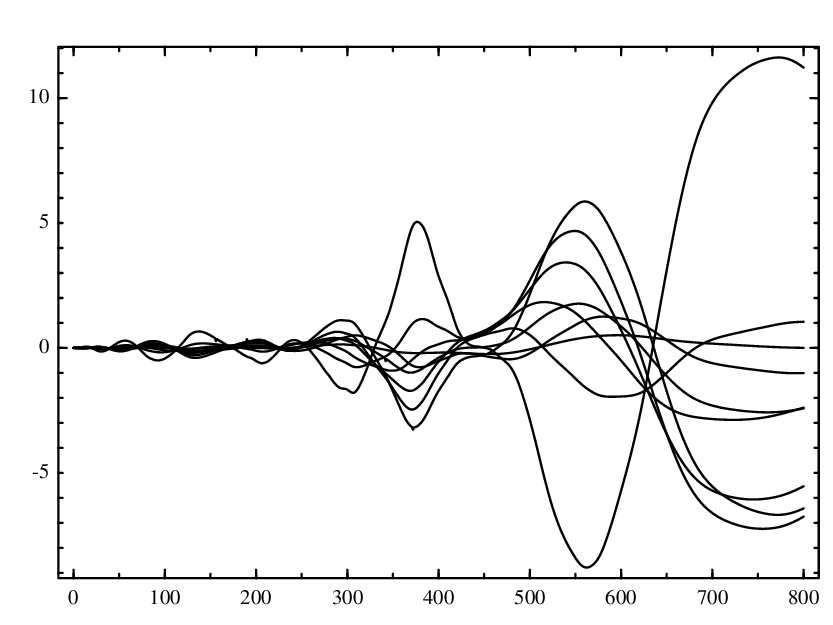}
\caption{$n=9$, left the state trajectory, all computations performed using 32 bits precision $\epsilon = 1.1921e-7$. Right, the error with 64 bits precision arithmetic $\epsilon = 2.2204e-16$}
\end{figure}

\begin{figure}
\includegraphics[width=6cm]{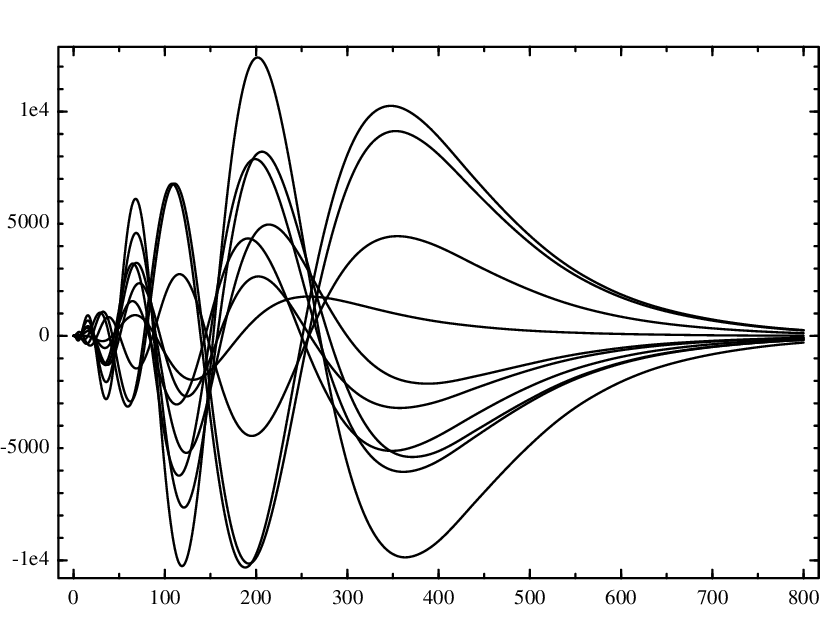}
\includegraphics[width=6cm]{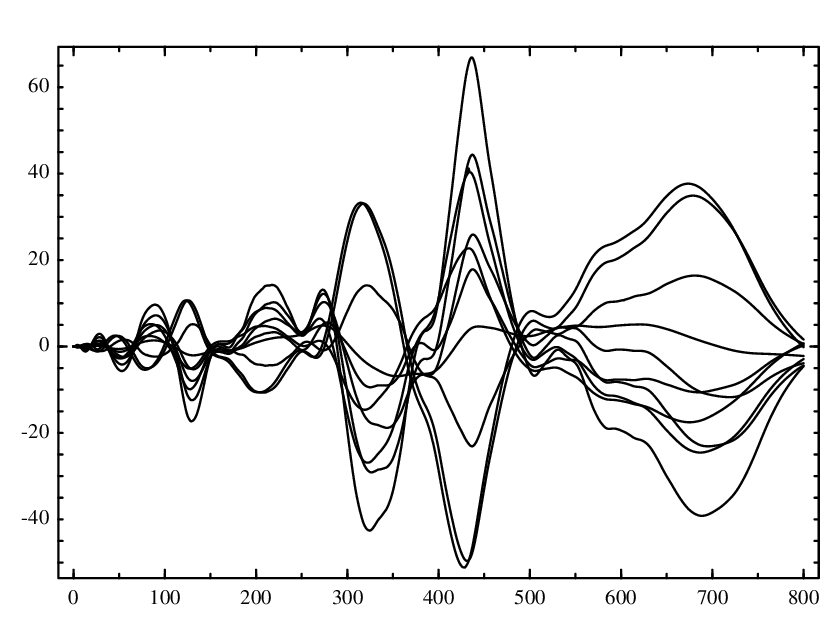}
\caption{$n=10$, left the state trajectory, all computations performed using 32 bits precision $\epsilon = 1.1921e-7$. Right, the error with 64 bits precision arithmetic $\epsilon = 2.2204e-16$}
\end{figure}

\begin{figure}
\includegraphics[width=6cm]{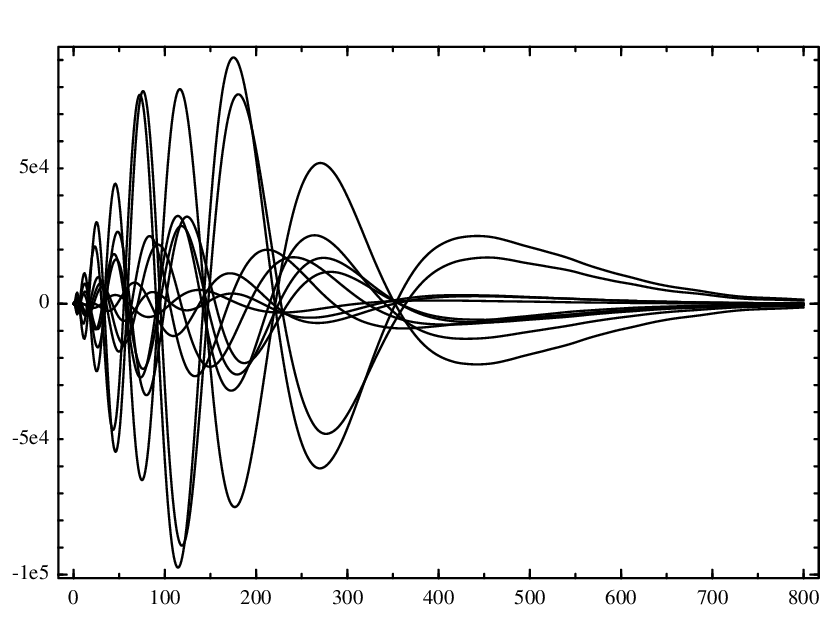}
\includegraphics[width=6cm]{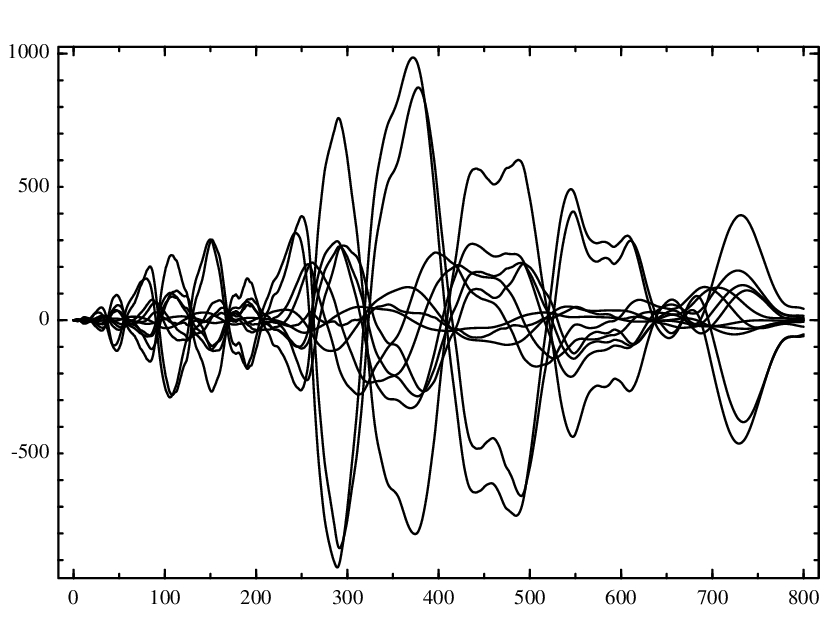}
\caption{$n=11$, left the state trajectory, all computations performed using 32 bits precision $\epsilon = 1.1921e-7$. Right, the error with 64 bits precision arithmetic $\epsilon = 2.2204e-16$}
\end{figure}

\begin{figure}
\includegraphics[width=6cm]{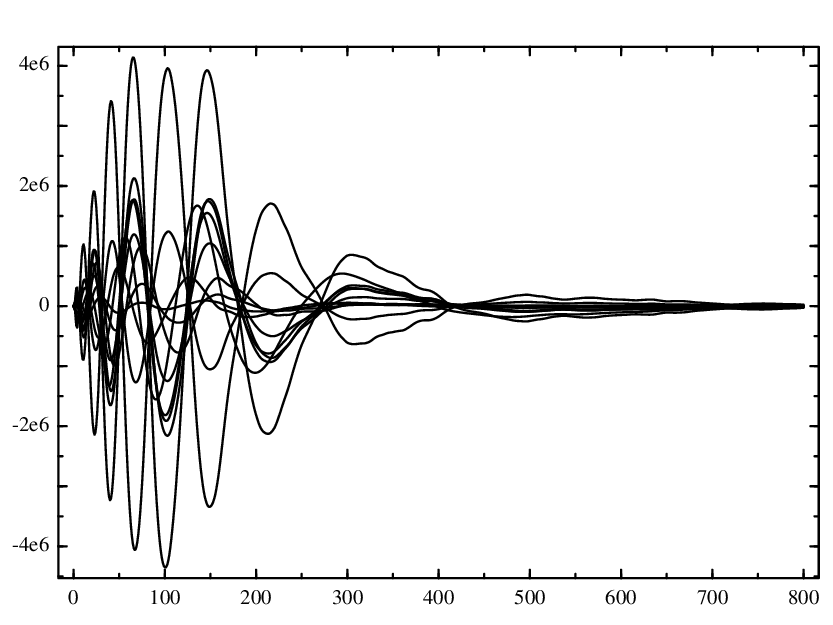}
\includegraphics[width=6cm]{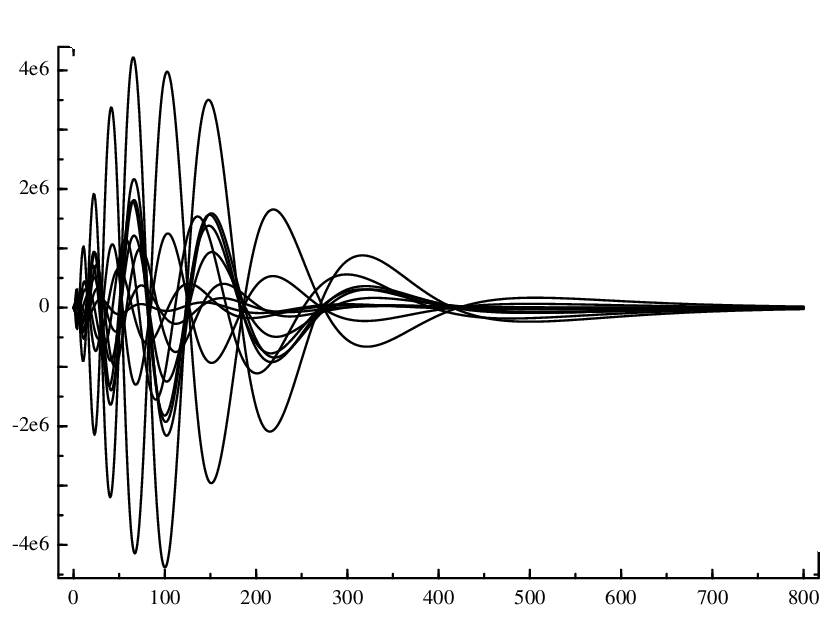}
\caption{$n=12$, left the state trajectory, all computations performed using 32 bits precision $\epsilon = 1.1921e-7$. Right, the state trajectory with 64 bits precision arithmetic $\epsilon = 2.2204e-16$}
\end{figure}

\begin{figure}
\includegraphics[width=6cm]{dim12_x32.eps}
\includegraphics[width=6cm]{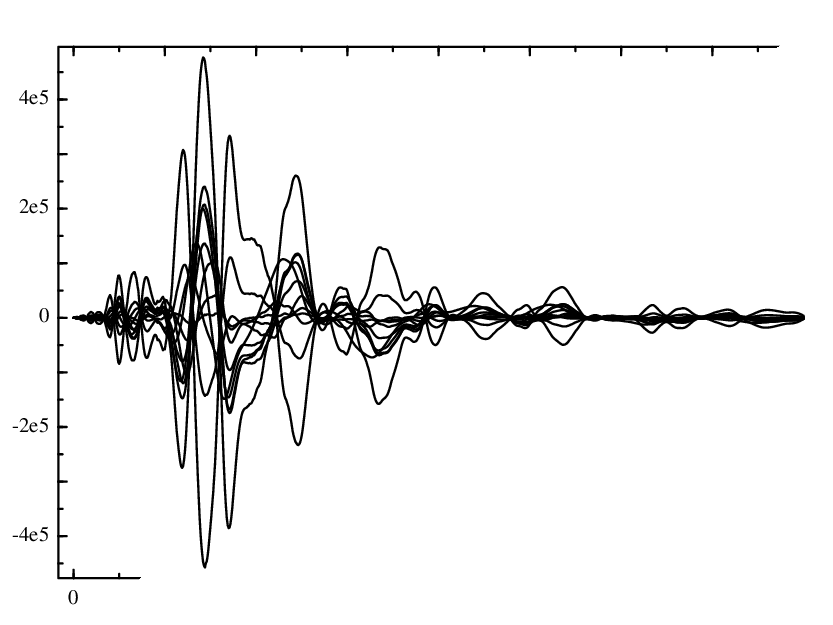}
\caption{$n=12$, left the state trajectory, all computations performed using 32 bits precision $\epsilon = 1.1921e-7$. Right, the error with 64 bits precision arithmetic $\epsilon = 2.2204e-16$}
\end{figure}

\begin{figure}
\includegraphics[width=6cm]{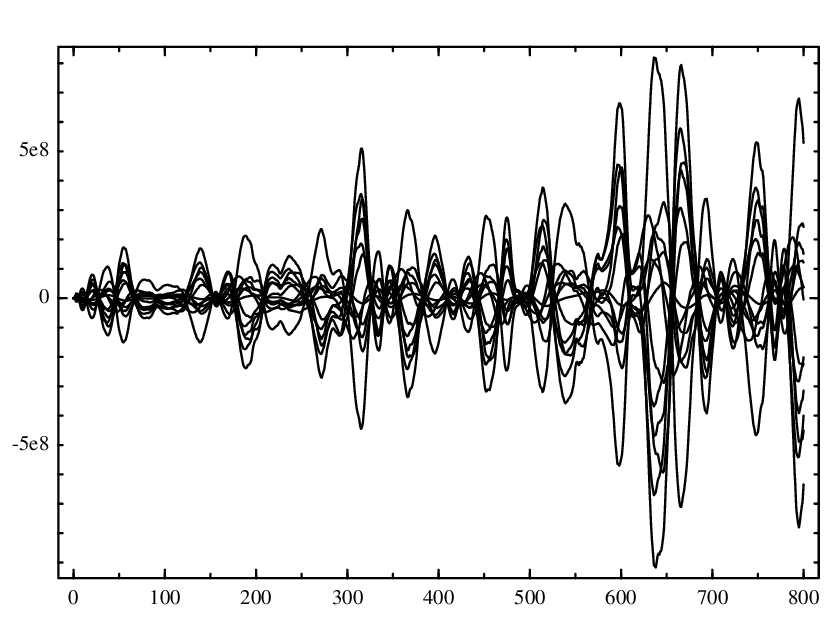}
\includegraphics[width=6cm]{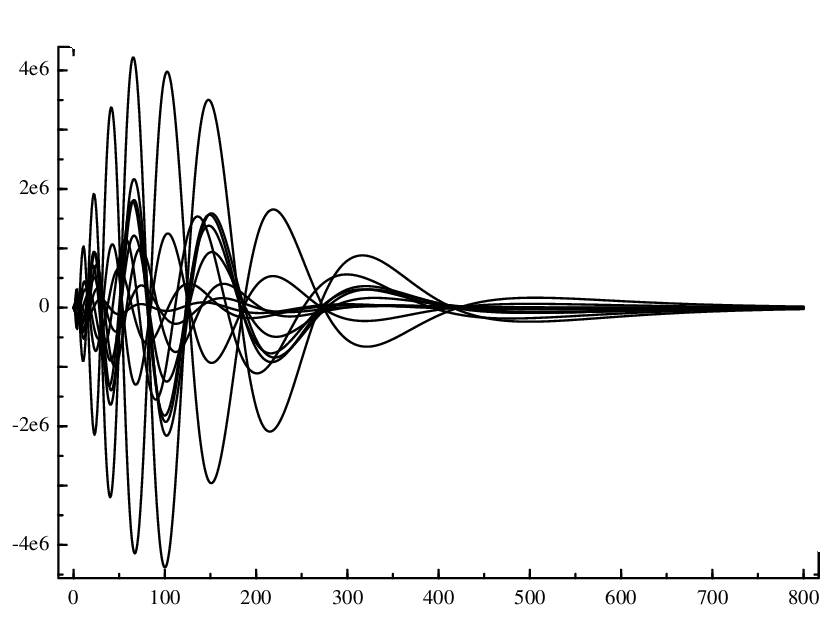}
\caption{$n=13$, (32 bits precision $\epsilon = 1.1921e-7$); left the state trajectory; right the state trajectory with 64 bits precision arithmetic $\epsilon = 2.2204e-16$}
\end{figure}

\section{Appendix: Computer Codes}

\begin{itemize}
\item The algebroid method
\item An example for the algebroid method
\item An example by intersection of hyperplanes
\item The ring operations method
\item The Miminis-Paige
\item The Varga algorithm
\item buildOp
\item fbkComp
\item fbkCompExtract
\item runMu3rk
\end{itemize}

\subsection{Example for the sliding intersection of hyperplances method}
\label{exampleAlgebroid}

\begin{verbatim}
%%%
%
% illustrates the sliding intersection of the affine hyperplanes of gains
%
A = [1 3 5; 7 13 17; 1 1 1];
B = [1; 1; 1];

% let us place the eigenvalues at -1, -2, and - 3
l1 = -1;
l2 = -2;
l3 = -3;
E = eye(3);
k11 = 1/(B(1))*(A(1,:)-l1*E(1,:));
k12 = 1/(B(2))*(A(2,:)-l1*E(2,:));
k13 = 1/(B(3))*(A(3,:)-l1*E(3,:));

k21 = 1/(B(1))*(A(1,:)-l2*E(1,:));
k22 = 1/(B(2))*(A(2,:)-l2*E(2,:));
k23 = 1/(B(3))*(A(3,:)-l2*E(3,:));

k31 = 1/(B(1))*(A(1,:)-l3*E(1,:));
k32 = 1/(B(2))*(A(2,:)-l3*E(2,:));
k33 = 1/(B(3))*(A(3,:)-l3*E(3,:));

% compute the vectors orthogonal to the hyperplanes 
n1 = inv([k11; k12; k13])*[1;1;1];
n2 = inv([k21; k22; k23])*[1;1;1];
n3 = inv([k31; k32; k33])*[1;1;1];

% compute the various normals (not normalized normals)
n10 = n1;
n20 = n2;
n30 = n3;

n21 = (eye(3) - n10*n10'*1/(n10'*n10))*n20;
n31 = (eye(3) - n10*n10'*1/(n10'*n10))*n30;

n32 = (eye(3) - n21*n20'*1/(n20'*n21))*n31;

G1 = n21*n20'*1/(n20'*n21);
G2 = n32*n30'*1/(n30'*n32);

gam1 = k11;
gam2 = gam1 - (gam1-k21)*G1';
gam3 = gam2 - (gam2-k31)*G2';

K = gam3
\end{verbatim}

\subsection{Example with the determinental intersection of hyperplanes}
\label{exampleIntersection}
The same numerical example as the previous section but relying on
intersection of hyperplanes using determinantal properties (for the determinantal properties used,  see for example no. 251 p.221 and no. 257 p. 226 of \cite{Dostor}).

\begin{verbatim}
A = [1 3 5; 7 13 17; 1 1 1];
B = [1; 1; 1];

% let us place the eigenvalues at -1, -2, -3
l1 = -1;
l2 = -2;
l3 = -3;
E = eye(3);

k11 = 1/(B(1))*(A(1,:)-l1*E(1,:));
k12 = 1/(B(2))*(A(2,:)-l1*E(2,:));
k13 = 1/(B(3))*(A(3,:)-l1*E(3,:));

k21 = 1/(B(1))*(A(1,:)-l2*E(1,:));
k22 = 1/(B(2))*(A(2,:)-l2*E(2,:));
k23 = 1/(B(3))*(A(3,:)-l2*E(3,:));

k31 = 1/(B(1))*(A(1,:)-l3*E(1,:));
k32 = 1/(B(2))*(A(2,:)-l3*E(2,:));
k33 = 1/(B(3))*(A(3,:)-l3*E(3,:));

%%
% compute the equations of the hyperplanes
% D + A x + B y + C z = 0
%
D1 = det([k11; k12; k13]);
ABC1 = [k11-k12; k11-k13];
A1 = -det(ABC1(:,[2,3]));
B1 = det(ABC1(:,[1,3]));
C1 = -det(ABC1(:,[1,2]));

D2 = det([k21; k22; k23]);
ABC2 = [k21-k22; k21-k23];
A2 = -det(ABC2(:,[2,3]));
B2 = det(ABC2(:,[1,3]));
C2 = -det(ABC2(:,[1,2]));

D3 = det([k31; k32; k33]);
ABC3 = [k31-k32; k31-k33];
A3 = -det(ABC3(:,[2,3]));
B3 = det(ABC3(:,[1,3]));
C3 = -det(ABC3(:,[1,2]));

%%
% compute the intersection of the three hyperplanes 
% this is [xx, yy, zz]
%
xx=-det([D1 D2 D3; B1 B2 B3; C1 C2 C3]);
yy=-det([A1 A2 A3; D1 D2 D3; C1 C2 C3]);
zz=-det([A1 A2 A3; B1 B2 B3; D1 D2 D3]);

dd = det([A1 A2 A3; B1 B2 B3; C1 C2 C3]);

K = 1/dd*[xx yy zz];
\end{verbatim}

\subsection{The first algebroid algorithm}

The algebroid method is described in Section \ref{algebroidDescription}.
To summarize, the quotients are taken in such a way that the representatives are in one-to-one correspondence with the points of the hyperplane associated with the current eigenvalue to be fixed. In the descending phase of the algorithm, each eigenvalue is fixed one by one, one in each quotient, the others remaining undefined. The size of the quotients decreases with each step. The final gain is obtained in an ascending phase.  This phase starts with the smallest quotient and works its way up to the original dimension. This construction phase combines the prior image of the quotient gain with the gain that initially placed the eigenvalue before taking the quotient, thus fixing all the undefined eigenvalues in each quotient at their correct values.

\label{algebroidAlgorithm}
\begin{verbatim}
%%%%%%%%%%%%%%%%%%%%%%%%%%%%%%%%%%%%%
%
function [K] = algebroid(A,B,VP)
%function [K] = algebroid(A,B,VP)
%
% algebroid based pole placement 
% It does not require that the system be brought
% to upper Hessenberg form prior to setting
% the eigenvalues.
%

n = length(B);
Ab = A;
Bb = B;
ani = zeros(n*(n-1),n);
kos = zeros(n,n);

for i=1:n-1
[qa,ra]=qr(Bb);
[qs,rs]=qr((qa(2:end,:)*(Ab-VP(i)*eye(n-i+1)))');
koh = qs'(end,:);
ko = Bb'/(Bb'*Bb)*(Ab - VP(i)*eye(n-i+1))*koh'*koh;
anb = qs'(1:end-1,:);
ani((i-1)*n+1: (i-1)*n+n-i, 1: n-i+1) = anb;
kos(i,1:n-i+1)=ko;
Ab = anb*Ab*anb';
Bb = anb*Bb;
end;

K = (Ab-VP(n))/Bb;
for i=n-1:-1:1
anb = ani((i-1)*n+1:(i-1)*n+n-i,1:n-i+1);
ko = kos(i,1:n-i+1);
K = ko + K*anb;
end;

K=-K; %convention A+BK
\end{verbatim}

We now present the version that replaces one of the qr decompositions (in fact the computationally more gready one) with the solution of a system of linear equations. This is quite useful since this displays the numerically critical, potentially ill-conditioned step in the entire process of eigenvalue assignment. This step appears during application of the classical Miminis-Paige algorithm by a subtle combination of the first step decomposition and its resulting effect in the qr step of the pole placement algorithm. (Recall that the first step is the step of reduction to Hessenberg control form using the staircase algorithm.)

\begin{verbatim}

%%%%%%%%%%%%%%%%%%%%%%%%%%%%%%%%%%%%%
%
function [K] = algebroid2(A,B,VP)
%function [K] = algebroid2(A,B,VP)
%
% algebroid based pole placement
% It does not require that the system be brought to upper
% Hessenberg form prior to setting the eigenvalues.
% Contrary to 'algebroid', this requires solving
% a linear system of equations. This replaces one of the
% two qr decompositions.
%

n = length(B);
Ab = A;
Bb = B;

ani = zeros(n*(n-1),n);
kos = zeros(n,n);

for i=1:n-1
Bbt = (Ab - VP(i)*eye(n-i+1))\Bb;
ko = Bbt/(Bbt'*Bbt);
[qb,rb]=qr(Bbt);
anb = qb(2:end,:);
ani((i-1)*n+1: (i-1)*n+n-i, 1: n-i+1) = anb;
kos(i,1:n-i+1)=ko;
Ab = anb*Ab*anb';
Bb = anb*Bb;
end;

K = (Ab-VP(n))/Bb;
for i=n-1:-1:1
anb = ani((i-1)*n+1:(i-1)*n+n-i,1:n-i+1);
ko = kos(i,1:n-i+1);
K = ko + K*anb;
end;

K=-K; %convention A+BK
\end{verbatim}

\subsection{The second algebroid algorithm}

This algorithm uses the coefficients of the characteristic polynomial rather thanthe explicit eingenvalues. It proceeds in two sweeps first the computation of the anchors (termed the'\texttt{oo}' projection operators) with the function '\texttt{buildOp}' and then it either extracts the feedback gain with '\texttt{fbkCompExtract}' or one can use a specficic feedback function called '\texttt{fbkComp}'. 

\subsubsection{Code for building the operators}
\label{buildOpCode}

The modified Ackermann's-based-formula algorithm constructs successive projection operators and quotients operators (nested $A$ operators) that are stored in arrays.  The following code provides the details of the computation and the storage.

\begin{verbatim}
%%%%%%%%%%%%%%%%%%%%%%%%%%%%%%%%%%%
%
function [At,Bt,Pt,Ao,Po,Oo] = buildOp(A,B)
%function [At,Bt,Pt,Ao,Po,Oo] = buildOp(A,B)

% builds the projection operators and the multiplication operators

n = length(B);

Ao = zeros(n*(n-1),n);
Po = zeros(n*(n-1),n);
Oo = zeros(n*(n-1),n);

Bt = B; % contains the current image of B
At = A;
Pt = eye(size(A));

for i=1:n-1

if (1==0) %oblique
	[q,r]=qr(At*B);
	if length(Bt)>=3
		nn = (q(2:end,:)*Bt)'*q(2:end,:);
		w = q(1:end-1,:);
		an = w-1/(nn*Bt)*w*Bt*nn;
	else
		an = q(1,:) - 1/(q(2,:)*Bt)*q(1,:)*Bt*q(2,:);
	end;
	Oo(((i-1)*n+1):((i-1)*n+(n-i)),1:n-i+1) = an;
	Ao(((i-1)*n+1):((i-1)*n+(n-i)),:) = an*At;

	% save weighting times annihilator
	Po(((i-1)*n+1):((i-1)*n+(n-i)),:) = an*Pt;

	% compute the image of B
	Bt = an*At*B;

	% compute the next nested A operator
	At = an*At*A;
	Pt = an*Pt;
else % orthogonal
% compute annihilator of B
	[Q,R] = qr(Bt);
	an = Q(2:end,:);
% find scaling by svd of operator times annihilator times A
% [u,d,v]=svd(an*At*B); %for Hessenberg Miminis Paige like 
	[u,d,v]=svd(an*At); % Mullhaupt version
	if (i<n-1)
		%w = diag(1./diag(d))*u';
		w = u';
		%w = eye(length(diag(d)));
	else
		%w = 1/d(1)*u';
		w = u';
		%w = 1;
	end;

	Oo(((i-1)*n+1):((i-1)*n+(n-i)),1:n-i+1) = w*an;
	Ao(((i-1)*n+1):((i-1)*n+(n-i)),:) = w*an*At;

	% save weighting times annihilator
	Po(((i-1)*n+1):((i-1)*n+(n-i)),:) = w*an*Pt;

	% compute the image of B
	Bt = w*an*At*B;

	%compute the next nested A operator
	At = w*an*At*A;
	Pt = w*an*Pt;

end;
end;

\end{verbatim}

\subsubsection{Code fbkComp}
\label{fbkCompCode}

Instead of using a single row of gain $u = - k x$ a feedback function $u = - k(x)$ is used that evaluates the value of $u$ based on the current $x$ using the nested quotients by performing more additions and multiplications. This improves the numerical accuracy. Below  is the function to do such computations.

\begin{verbatim}
%%%%%%%%%%%%%%%%%%%%%%%%%%%%%%%
%
function u = fbkComp(Ao,Oo,A,B,P,x)
%function u = fbkComp(Ao,Oo,A,B,P,x)

% poles are in P

pp = poly(P);
pp = pp(end:-1:1); % put in reverse order

ut = pp(1)*x; % cumulative contribution of each image to the input

n = length(B);

for i=1:n-1
Ot = Oo(((i-1)*n+1):((i-1)*n+(n-i)),1:n-i+1);
At = Ao(((i-1)*n+1):((i-1)*n+(n-i)),:);
ut = At*x*pp(i+1) + Ot*ut;
end;
ut = ut+At*A*x;
u = -1/(At*B)*ut;
\end{verbatim}

\subsubsection{Code fbkCompExtract for extracting the $K$ vector} 
\label{fbkCompExtractCode}

It is possible to compute the $K$ vector from the nested $A_i$ operators and projection operators $O_i$. The following function returns the value of the row vecor of gains $K$ such that $u = - K x$ (a classical scalar product).

\begin{verbatim}
function K = fbkCompExtract(Ao,Oo,A,B,P)
%function K = fbkCompExtract(Ao,Oo,A,B,P)

% poles are in P

n = length(B);

pp = poly(P);
pp = pp(end:-1:1); % put in reverse order

Kt = pp(1)*eye(n); % cumulative contribution to the final gain 

for i=1:n-1
Ot = Oo(((i-1)*n+1):((i-1)*n+(n-i)),1:n-i+1);
At = Ao(((i-1)*n+1):((i-1)*n+(n-i)),:);
Kt = At*pp(i+1) + Ot*Kt;
end;
Kt = Kt+At*A;
K = -1/(At*B)*Kt;


\end{verbatim}

\subsubsection{Code for runMu3rk}

This code can be used to test the placement of eigenvalues by integrating a differential equations and testing its trajecteries.

\begin{verbatim}
%%%%%%%%%%%%%%%%%%%%%%%%%%%%%%%%%%%%%
%
function [t,x,A,B,Ao,Oo,xm] = runMu3rk(n,T,h);
%function [t,x,A,B,Ao,Oo,xm] = runMu3rk(n,T,h);


[A,B] = buildEx(n);
[At,Bt,Pt,Ao,Po,Oo] = buildOp(A,B);

K = fbkCompExtract(Ao,Oo,A,B,-(1:n)*0.01); %for proposed version
%K = miminis(A,B,-(1:n)*0.01); %for Miminis version
%
% Remark: uncomment one of the above two lines
%         to activate the desired version of pole placement

t=0:h:T;
t=t';
N = length(t);
x = zeros(N,n);
x(1,:) = 1:n;
for i=1:N-1
k1 = simMu3(t,x(i,:)',A,B,Ao,Oo,-(1:n)*0.01,K,1)';
k2 = simMu3(t,x(i,:)' + k1'*h/2,A,B,Ao,Oo,-(1:n)*0.01,K,1)';
k3 = simMu3(t,x(i,:)' + k2'*h/2,A,B,Ao,Oo,-(1:n)*0.01,K,1)';
k4 = simMu3(t,x(i,:)' + k3'*h,A,B,Ao,Oo,-(1:n)*0.01,K,1)';
x(i+1,:) = 1/6*(k1 + 2*k2+ 2*k3 + k4)*h + x(i,:);
end;
clf

subplot(311)
plot(t',x',{LineWidth=0.3});

xm = zeros(N,n);
xm(1,:) = 1:n;
for i=1:N-1
k1 = simMu3(t,xm(i,:)',A,B,Ao,Oo,-(1:n)*0.01,K,0)';
k2 = simMu3(t,xm(i,:)' + k1'*h/2,A,B,Ao,Oo,-(1:n)*0.01,K,0)';
k3 = simMu3(t,xm(i,:)' + k2'*h/2,A,B,Ao,Oo,-(1:n)*0.01,K,0)';
k4 = simMu3(t,xm(i,:)' + k3'*h,A,B,Ao,Oo,-(1:n)*0.01,K,0)';
xm(i+1,:) = 1/6*(k1 + 2*k2+ 2*k3 + k4)*h + xm(i,:);
end;

subplot(312)
plot(t',xm',{LineWidth=0.3});

subplot(313)
plot(t',x'-xm',{LineWidth=0.3});

function xd = simMu3(t,x,A,B,Ao,Oo,P,K,toggle)
%function xd = simMu3(t,x,A,B,Ao,Oo,P,K,toggle)

if toggle==1
u = fbkComp(Ao,Oo,A,B,P,x);
else
u = K*x;
end;

xd = A*x + B*u;
\end{verbatim}

\subsection{The second algebroid algorithm, the ring operations method}

A method that uses mainly ring operations is now provided in Mathematica (because of the need to have infinite digits integers). It requires the GCD of two elements of the ring (in implementing the NullSpace method of Mathematica). Since the implementation of NullSpace is hidden, we managed to reverse engineer the method and NullSpace2 is a function leading to similar results as NullSpace. It requires the computation of the GCD between two elements of the ring. Hence the method requires some division process. The GCD limits the growth of the elements of the ring. So formally, the method proposed works for an Euclidean ring. 

The gain vector $K$ is synthesized step by step while computing the quotients rather than proceeding in a two sweep method.

\subsubsection{The main ring operator algorithm, Mathematica version}

\begin{verbatim}
(* Mathematica script *)
PlaceAlg[A_, B_, P_] := Module[{Ab, Bb, anb, As, Bs, ans, PP, KK},
   Ab = IdentityMatrix[Length[A]];
   Bb = B;
   As = {Ab};
   Bs = {Bb};
   ans = {};
   PP = P;
   KK = IdentityMatrix[Length[A]]*First[PP];
   PP = Rest[PP];
   Do[Block[{},
      anb = NullSpace[Transpose[Bb]];
      As = Append[As, Ab];
      Bs = Append[Bs, Bb];
      ans = Append[ans, anb]; 
      Bb = anb . Ab . A . B;
      Ab = anb . Ab . A;
      KK = First[PP] Ab + anb . KK;
      PP = Rest[PP]];, Length[A] - 1];
   KK = KK + Ab . A;
   {(Ab . B // Flatten // Part[#, 1] &), KK}
   ];
\end{verbatim}

\subsubsection{The NullSpace2 function}

Instead of the 'NullSpace' function appearing just afte the 'Do' statement, the following function performs similarly than NullSpace (concerning the integer explosion). If the author is not misled the 'NullSpace' differs from 'NullSpace2' by systematically taking the same pivotal element and scanning the rest of the elements sequentially. 'NullSpace2' takes sequentially pairwise elements. The important step is to take the greatest common divisor between the chosen elements. Hence it introduces divisions, both when dividing the elements by the GCD and when computing the GCD --- if one takes the classical Euclidean algorithm using division modulo arithmetic. There exists GCD algorithms that do not use divisions only substractions. They are not as efficient as division algorithms. 
\begin{verbatim}
NullSpace2[BB_] := Module[{B, annil, NN, li, j},
   B = BB[[1]];
   NN = Length[B];
   li = Table[0, {j, 1, NN}];
   annil = {};
   For[j = 1, j < NN, j++, Module[{lit},
      lit = li;
      If[B[[j]] == 0,
       Module[{},
         lit[[j]] = 1; 
         ];,
       Module[{k},
         k = j + 1;
         While[(B[[k]] == 0 ) && (k < NN), k = k + 1;];
         If[B[[k]] != 0, 
          Module[{gcd}, gcd = GCD[B[[k]], B[[j]]]; 
            lit[[k]] = -B[[j]]/gcd;
            lit[[j]] = B[[k]]/gcd;];,
          lit[[k]] = 1;];
         ];];
       annil = Append[annil, lit];];];
   annil];
\end{verbatim}
\begin{remark}
We do not claim that 'NullSpace2' is a replacement for 'NullSpace' for a general matrix. It is a replacement when dealing with a single row vector.
\end{remark}
To be complete, here are two versions of the GCD algorithm to replace 'GCD':
\begin{verbatim}
(* substraction version *)
GCDsub[a_, b_] := If[a==b, b, If[a>b,GCDsub[a-b,b],GCDsub[a,b-a]]];

(* division modulo version *)
GCDmod[a_, b_] := If[b==0, a, GCDmod[b, Mod[a, b]]];
\end{verbatim}

\subsubsection{A Gambit Scheme implementation}

If one does not want to rely on commercial software like Mathematica, a Scheme version is now proposed. It takes advantage of the the large integers implementation. It has been checked with the Gambit scheme implementation. Warning: there might be issues while copying and pasting from the .pdf directly. It is advised to type each line seperatly by hand. 
\small
\begin{verbatim}
;; false pseudo-random generator (cheap but easy to cross-code)
(define seed 341)
(define falserandom 
  (lambda (m) (set! seed 
               (modulo (* seed 564485031) 1234465738)) 
              (quotient (modulo seed (* m 123)) 123)))

;; same as before but negative numbers as well
(define centeredfalserandom (lambda (m) (- (falserandom (* 2 m)) m)))

;; (gcdMod a b) computes the greates common divisor using Euclid's algo
(define gcdMod (lambda (a b) (if (= b 0) a (gcdMod b (modulo a b)))))

;; (gcdMod a b) computes the greatest common divisor using only substractions
(define gcdSub (lambda (a b) (if (= a b) b (if (> a b) (gcdSub (- a b) b) 
(gcdSub a (- b a))))))

;; (genCons n val) generates a vector of length n of constant values val
(define genCons (lambda (n val) 
   (letrec ((x (lambda (n li) 
                (if (or (zero? n) (< n 0)) li (x (- n 1) (cons val li)))))) 
            (x n '()))))

;; (genB n) generates a vector of ones of length n. This is the
;; the B vector of a few examples
(define genB (lambda (n) (genCons n 1)))

;; 
(define genA (lambda (n) (map (lambda (x) (reverse x)) (reverse (eyep n 2)))))


;; (genZ n) creates a vector of zeros of length n.
(define genZ (lambda (n) (genCons n 0)))

;; (zeros n m) creates an array of zeros n x m
(define zeros (lambda (n m) 
   (letrec ((x (lambda (n m li)  (if (or (zero? n) (< n 0)) li 
                                  (x (- n 1) m (cons (genCons m 0) li)))))) 
            (x n m '()))))

(define A0 '((1 3 5) (7 13 17) (1 1 1)))

(define B0 '(1 1 1))

;; scalar times vector
(define scv (lambda (s v) (map (lambda (x) (* s x)) v)))

;; scalar times cracovian or matrix
(define scm (lambda (s cm) (map (lambda (x) (scv s x)) cm)))
 
;; scalar product between vector v1 and vector v2
(define scal (lambda (v1 v2) 
  (letrec ((x (lambda (v1 v2 sc) 
              (if (null? v1) sc (x (cdr v1) (cdr v2) 
              (+ sc (* (car v1) (car v2)))))))) 
           (x v1 v2 0))))

;; scal2 is slower in execution time
(define scal2 (lambda (v1 v2)
(apply + (map (lambda (x y) (* x y)) v1 v2))))

;; product of matrix A times vector B
;; or row multiplication of each entry of cracovian A times vector B
(define prod (lambda (A B) (map (lambda (x) (scal x B)) A)))

;; multiplication of cracovians
;; the final result is the transposed value of the standard cracovian
;; i.e. row product, row cracovian, instead of column cracovian
(define prodc (lambda (a b) (map (lambda (x) (prod a x)) b)))

;; transpose of a cracovian (or matrix) is obtained by postmultiplying
;; by the elementary matrix
(define tT (lambda (m) (prodc m (eye (length (car m))))))

;; addition of vectors
(define addv (lambda (x y) (map (lambda (x y) (+ x y)) x y)))

;; addition of cracovians or matrices
(define addcm (lambda (x y) (map addv x y)))

;; append to a list li the last element el
(define append (lambda (li el) (if (null? li) (list el) (cons (car li) (append (cdr li) el)))))

;; multiply the monomial mono to the polynomial poly
(define multmp (lambda (mono poly) 
(let 
  ((p1 (append (scv (car mono) poly) 0))
   (p2 (cons 0 (scv (cadr mono) poly))))
  (map (lambda (x y) (+ x y)) p1 p2))))

;; poly creates the polynomial whose roots are the elements of li
(define poly (lambda (li) (letrec ((poly (lambda (lin lout) 
(if (null? lin) lout (poly (cdr lin) (multmp (list 1 (- (car lin))) lout))))))
(poly li '(1)))))

;; suppresses the first column of matrix M, returns the remaining matrix
(define rest (lambda (M) (map (lambda (x) (cdr x)) M)))

;; gets the first column of matrix M and return it as a vector
(define first (lambda (M) (map (lambda (x) (car x)) M)))

;; identityMatrix
(define eye (lambda (n)  
  (letrec ((x (lambda (m mat) 
               (if (= m n) mat 
                           (x (+ m 1) (cons (cons 1 (genZ m)) 
                                       (map (lambda (x) (cons 0 x)) mat)))))))
           (x 1 '((1))))))

;; same is identity matrix but with a particular value v instead of 1
;; on the main diagonal 
(define eyev (lambda (n v)  (letrec ((x (lambda (m mat) (if (= m n) mat (x (+ m 1)
                          (cons (cons v (genZ m)) 
                                (map (lambda (x) (cons 0 x)) mat)))))))
                                 (x 1 (list (list v))))))

(define eyep (lambda (n p)  (letrec ((x (lambda (m mat v) (if (= m n) mat (x (+ m 1)
                          (cons (cons (expt p v)  (genZ m)) 
                                (map (lambda (x) (cons 0 x)) mat)) (+ v 1))))))
                                 (x 1 (list (list 1)) 1))))

;; nullspace
(define nullspace (lambda (BB) 
(let ((dupliczero (lambda (x) (cons (map (lambda (x) 0) (car x)) x)))
      (augment (lambda (x) (map (lambda (x) (cons 0 x)) x)))
      (insert (lambda (a b x) (cons (cons (- b) (cons a (cddar x))) (cdr x)))))
      (letrec ((nullsp (lambda (anni b) 
                            (if (null? (cdr b)) 
                              anni 
                              (nullsp (let ((x (cadr b)) 
                                             (y (car b)) 
                                          (gcd (gcdMod (cadr b) (car b)))) 
                                          (insert (quotient x gcd) (quotient y gcd) 
                                                (augment (dupliczero anni))))
                                      (cdr b))))))
               (map (lambda (x) (reverse x)) 
                   (nullsp (let ((a (car BB)) 
                                 (b (- (cadr BB))))
                              (let ((gcd (gcdMod a b)))
                                 (list (list (quotient a gcd) (quotient b gcd))))) (cdr BB)))))))

;; BB with random numbers
(define BB (map (lambda (x) (falserandom 300)) (genZ 15))) 

(define placeMu (lambda (A0 B0 PP) 
(let ((AB (prod A0 B0)) (A0t (tT A0)))
  (letrec ((placeMu (lambda (Bb Ab KK PP) 
                   (if (null? (cdr PP)) 
		     (list (prod Ab B0) (addcm KK (prodc A0t Ab)))
                     (let* ((anb (nullspace Bb))
                            (Ab2 (prodc (prodc Ab A0t) anb)))
		(begin ;(display "Ab = ") (display Ab) (display "\n")
		       ;(display "Ab2 = ") (display Ab2) (display "\n")
                       ;(display "anb = ") (display anb) (display "\n")
                       ;(display "Bb = ") (display Bb) (display "\n")
                     (placeMu (prod anb (prod Ab AB)) 
                              Ab2 
			      (addcm (scm (car PP) Ab2) (prodc (tT KK) anb))
                              (cdr PP))))))))
           (placeMu B0 
                    (eye (length A0)) 
                    (eyev (length A0) (car PP)) 
                    (cdr PP))))
))

(define simpli (lambda (KK) 
  (let* ((coef (caar KK))
         (vect (caadr KK))
         (gcdall (apply gcd (cons coef vect))))
    (list (list (quotient coef gcdall)) 
        (list (map (lambda (x) (quotient x gcdall)) vect))))))

(define ratio (lambda (KK)
     (car (scm (/ 1 (caar KK)) (cadr KK)))))

(define progress (lambda (n li)
 (if (zero? n) li (progress (- n 1) (cons n li)))))

(define genP (lambda (n)
   (reverse (poly (map (lambda (x) (- x)) (progress n '()))))))

(define KK (placeMu A0 B0 '(6 11 6 1)))

(define ratKK (ratio KK))

(define A5 (genA 5))
(define B5 (genB 5))
(define P5 (reverse (poly '(-1 -2 -3 -4 -5))))

(define A7 (genA 7))
(define B7 (genB 7))
(define P7 (reverse (poly '(-1 -2 -3 -4 -5 -6 -7))))

;(define KK5 (placeMu A5 B5 P5))
;(define ratKK5 (scm (/ 1 (caar KK5)) (cadr KK5)))

;(define KK7 (placeMu A7 B7 P7))
;(define ratKK7 (scm (/ 1 (caar KK7)) (cadr KK7)))

(display (ratio (simpli KK)))

;;(ratio (simpli (placeMu (genA 150) (genB 150) (genP 150))))
\end{verbatim}

\subsection{The Miminis-Paige algorithm}
\label{MiminisCode}
\begin{verbatim}

%%%%%%%%%%%%%%%%%%%%%%%%%%%%%%%%%%%%%
%
function [K,A1,B1,Kb,qc] = miminis(A,B,VP)
%function [K,A1,B1,Kb,qc] = miminis(A,B,VP)

n = length(B);

% simultaneous preparation of B and upper Hessenberg A 
temp = [B,A];
qc = eye(n);
VP = VP(end:-1:1);
for i=1:301
[q,r]=qr(temp);
qc = qc*q;
temp = q'*temp*[1, zeros(1,n); zeros(n,1), q];
end;

Ai = (qc'*A*qc)(end:-1:1,end:-1:1);
Bi = (qc'*B)(end:-1:1);
A1=Ai;
B1=Bi;

qqi = zeros(n*(n-1),n);
pph = zeros(1,n);

% forward phase, nested Hessenberg modification
for i=1:n-1
[qi,ri]=qr(Ai'-VP(i)*eye(n-i+1));
pph(i) = ri(end,end)/Bi(end);
qqi((i-1)*n+1:(i-1)*n+(n-i+1),1:n-i+1)=qi;
if i<n-1
Ai = (qi'*Ai*qi)(1:end-1,1:end-1);
Bi = (qi'*Bi)(1:end-1);
else
Ai = (qi'*Ai*qi)(1,1);
Bi = (qi'*Bi)(1);
end;
end; %for
pph(n) = (Ai-VP(n))/Bi;
K=[pph(n)];

% backward gain-construction phase 
for i=n-1:-1:1
qi = qqi((i-1)*n+1:(i-1)*n+(n-i+1),1:n-i+1);
K=[K,pph(i)]*qi';
end;

Kb = K;
K = -K(end:-1:1)*qc';
\end{verbatim}

\subsection{The Varga algorithm}
\label{VargaAlgo}

Below is an implementation of the Varga algorithm in Matlab-like syntax (SysQuake).

\begin{verbatim}
function K = schurVarga(A,B,P)
% function K = schurVarga(A,B,P)
% 
% pole placement using the Varga method through schur decomposition
% A = control system matrix
% B = single column vector associated with input
% P = vector of eigenvalues to be placed
%
% K is the vector of state feedback gain such that A + B*K has P as
% eigenvalues
%
(U,T) = schur(A);
As = T;
Bs = U'*B;
ii = length(B);
n = ii;
hh = zeros(1,n);
Qs = eye(n);
while ii>1
ii
abs(As(ii,ii-1))<5
if abs(As(ii,ii-1))<5*eps % single real eigenvalue
  k1 = 1/Bs(end)*(P(ii) - As(end,end));
  kk = [zeros(1,n-1),k1];
  hh = hh + kk*Qs;
  As = As+Bs*kk
  % now bring down the top most unasigned diag element down to
  % lowest position using successive permutations of 2 along
  % the diagonal
  Qs = eye(n);
  for jj = 2:n
    (As, Q) = exchgIJ(As,jj,n);
    Qs = Q*Qs;
  end;
  Bs = Qs*Bs;
  ii = ii-1;
else % pair of complex conjugate eigenvalues
  A2 = As(ii-1:ii,ii-1:ii);
  B2 = Bs(ii-1:ii);
  (U2,B2,V2) = svd(B2);
  A2s = U2'*A2*U2;
  k1s = (P(ii-1) + P(ii) - A2s(1,1) - A2s(2,2))/B2s(1);
  k2s = (A2s(2,2)/A2s(2,1))*k1s + (A2s(1,1)*A2s(2,2)-A2s(1,2)*A2s(2,1) -...
	 P(ii-1)*P(ii))/(A2s(2,1)*B2s(1));
  k2 = V2*[k1s, k2s]*U2';
  kk = [kk,k2];
  hh = kk;
  ii = ii-2;
end;
end;

K = hh;

function (Ks,As,Bs,K) = schurRVarga(A,B,P)
% function (Ks,As,Bs,K) = schurRVarga(A,B,P)
% 
% Pure real eigenvalues Varga method, both the system A and desired
%  eignvalues are real.
% pole placement using the Varga method through schur decomposition
% A = control system matrix with purely real eigenvalues
% B = single column vector associated with input
% P = vector of eigenvalues to be placed they are all real eigenvalues.
%
% K is the vector of state feedback gain such that A + B*K has P as
% eigenvalues
%

(U,T) = schur(A);
As = T;
Bs = U'*B;
Bsd = Bs;
ii = length(B);
n = ii;
hh = zeros(1,n);
Qs = eye(n);
while ii>0
  k1 = 1/Bs(end)*(P(ii) - As(end,end));
  kk = [zeros(1,n-1),k1];
  hh = hh + kk*Qs;
  As = As+Bs*kk;
  % now bring down the top most unasigned diag element down to 
  % lowest position using successive permutations of 2 along
  % the diagonal
  if ii>=1
    for jj = 2:n
      (As, Q) = exchgIJ(As,jj,n);
      Qs = Q*Qs;
    end;
    Bs = Qs*Bsd;
  end;
  ii = ii-1;
end;
Ks = hh;
As = T;
Bs = U'*B;
eig(As+Bs*Ks)
K = Ks*U';

function (So, Q, iter) = exchg2it(Si)
% function (So, Q, iter) = exchg2it(Si)
%
% exchanges two elements of the diagonal of Si using
% an orthogonal transformation while leaving the Schur form intact
% So is the resulting output Schur form. 
% Si is the inital Schur form.
% Q contains the associated similarity transform
% A succession of Givens rotations are used.
%
P = [0,1;1,0];
Q = P;
H = P*Si*P;
iter = 0;
while abs(H(2,1)) > 0.5*eps
s = H(2,1)/sqrt(H(:,1)'*H(:,1));
c = H(1,1)/sqrt(H(:,1)'*H(:,1));
G = [c, s; -s, c];
H = G*H*G';
Q = G*Q;
iter = iter + 1;
end;
So = H;

function (So, Q) = exchg2(Si)
% function (So, Q) = exchg2(Si)
%
% exchanges two elements of the diagonal of Si. 
% using an orthogonal transformation while keeping the Schur structure.
% So is the Schur form after transformation.
% Si is the initial Schur form.
% Q contains the associated similarity transform.
% A single givens Rotation is used.
%
P = [0,1;1,0];
Q = P;
H = P*Si*P;
c = H(1,1)-H(2,2);
s = H(2,1);
den = sqrt(c^2+s^2);
c = c/den;
s = s/den;
G = [c, s; -s, c];
H = G*H*G';
Q = G*Q;
So = H;

function (So, Q) = exchgIJ(Si, ii, n)
% function (So, Q) = exchgIJ(Si, ii, n)
P = eye(n);
P(ii-1:ii, ii-1:ii) = [0,1;1,0];
Q = P;
H = P*Si*P;
c = H(ii-1, ii-1)-H(ii,ii);
s = H(ii,ii-1);
den = sqrt(c^2+s^2);
c = c/den; s = s/den;
G = eye(n);
G(ii-1:ii, ii-1:ii) = [c, s; -s, c];
H = G*H*G';
Q = G*Q;
So = H;

\end{verbatim}

\section{Conclusion}
 A framework for pole placement for single-input linear time-invariant systems is provided with a new variant of the pole placement algorithm (Miminis-Paige like algorithm without requiring an Hessenberg initial reduction). The framework is the use of quotient operations and choices of representatives of the equivalence class.

When the eigenvalues to be placed are real and distinct from each other and distinct from the eigenvalues of $A$, the framework gives a geometric intepretation of the gain vector as the vector corresponding to the intersection of $n-1$ skew affine (i.e. not passing through the origin) hyperplanes. Each affine hyperplane is associated with one of the eigvenvalues to be placed (supposed all to be distinct leading to non parallel independent hyperplanes). The fact that the eigenvalues are distinct from each other guarantees that the hyperplanes are independent and intersect at a single point giving the final gain vector.

Links and connections with classical literature algorithms of Miminis-Paige and Varga have been provided and discussed.

A factorization of Ackermann's formula can revalidate the approach to a certain extent to tackle ill-conditioned systems. This has been illustrated with a numerical ill-conditioned system. All codes are provided in full.

\section{Acknowledgment}

The author would like to extend his deepest gratitude to Dr. Yves Piguet for providing specific compiled LME and SysQuake 6.5 versions.

\section{Literature and comments}

\cite{CoxMoss}, \cite{CoxMoss2}, \cite{ArnoldDatta}, \cite{Boley}, \cite{BoleyGolub}, \cite{Varga},
\cite{MiminisPaige1}, \cite{MiminisPaige2}, \cite{Tsui}, \cite{Tsui1}, \cite{Hyman}, \cite{Patel},
\cite{Petkov}, \cite{Petkov2}, \cite{Paige}, \cite{PatelMisra}, \cite{KautskyNicholsVanDooren},
\cite{DattaElhayRam}, \cite{DattaSaad}, \cite{Datta}, \cite{BhattacharyyaDeSouza}, \cite{ArnoldDatta},
\cite{Ackermann}, \cite{Dostor}, \cite{Gantmacher1}, \cite{Gantmacher2}, \cite{Golub}, \cite{Wilkinson}

\bibliographystyle{plain}
\bibliography{sevasBib}

\end{document}